\documentclass[oneside,english]{amsart}
\usepackage[T1]{fontenc}
\usepackage[latin9]{luainputenc}
\usepackage{geometry}
\usepackage{array}
\usepackage{longtable}
\usepackage{float}
\usepackage{multirow}
\usepackage{amstext}
\usepackage{amsthm}
\usepackage{amssymb}
\usepackage{graphicx}
\usepackage{setspace}

\makeatletter

\providecommand{\tabularnewline}{\\}

\numberwithin{equation}{section}
\numberwithin{figure}{section}
\theoremstyle{plain}
\newtheorem{thm}{\protect\theoremname}
\theoremstyle{definition}
\newtheorem{defn}[thm]{\protect\definitionname}
\theoremstyle{plain}
\newtheorem{lem}[thm]{\protect\lemmaname}
\theoremstyle{definition}
\newtheorem{example}[thm]{\protect\examplename}
\theoremstyle{remark}
\newtheorem{rem}[thm]{\protect\remarkname}

\makeatother

\usepackage{babel}
\providecommand{\definitionname}{Definition}
\providecommand{\examplename}{Example}
\providecommand{\lemmaname}{Lemma}
\providecommand{\remarkname}{Remark}
\providecommand{\theoremname}{Theorem}

\begin{document}
\begin{singlespace}
\noindent \title[Centralizers of  nilpotent elements in basic classical  Lie superalgebras]{Centralizers of nilpotent elements in basic classical Lie superalgebras in good characteristic}
\end{singlespace}
\author{Leyu Han}
\maketitle
\noindent \begin{center}
\textit{School of Mathematics, University of Birmingham, }
\par\end{center}

\noindent \begin{center}
\textit{Birmingham, B15 2TT, UK }
\par\end{center}

\noindent \begin{center}
\textit{feish.ly@gmail.com}
\par\end{center}

\noindent \begin{center}
\textit{ORCID: 0000-0002-4170-4210}
\par\end{center}
\begin{abstract}
\noindent Let $\mathfrak{g}=\mathfrak{g}_{\bar{0}}\oplus\mathfrak{g}_{\bar{1}}$
be a basic classical Lie superalgebra over an algebraically closed
field $\mathbb{K}$ whose characteristic $p>0$ is a good prime for
$\mathfrak{g}$. Let $G_{\bar{0}}$ be the reductive algebraic group
over $\mathbb{K}$ such that $\mathrm{Lie}(G_{\bar{0}})=\mathfrak{g}_{\bar{0}}$.
Suppose $e\in\mathfrak{g}_{\bar{0}}$ is nilpotent. Write $\mathfrak{g}^{e}$
for the centralizer of $e$ in $\mathfrak{g}$ and $\mathfrak{z}(\mathfrak{g}^{e})$
for the centre of $\mathfrak{g}^{e}$. We calculate a basis for $\mathfrak{g}^{e}$
and $\mathfrak{z}(\mathfrak{g}^{e})$ by using associated cocharacters
{\normalsize{}$\tau:\mathbb{K}^{\times}\rightarrow G_{\bar{0}}$}
of $e$. In addition, we give the classification of $e$ which are
reachable, strongly reachable or satisfy the Panyushev property for
exceptional Lie superalgebras $D(2,1;\alpha)$, $G(3)$ and $F(4)$.
\end{abstract}



\section{Introduction\label{sec:Introduction}}

\noindent Let $\mathfrak{g}=\mathfrak{g}_{\bar{0}}\oplus\mathfrak{g}_{\bar{1}}$
be a basic classical Lie superalgebra over an algebraically closed
field $\mathbb{K}$ whose characteristic $p>0$ is a good prime for
$\mathfrak{g}$. Note that the definition for a good prime is a natural
extension of that for simple Lie algebras (see Definition \ref{def:good prime}).
Let $e\in\mathfrak{g}_{\bar{0}}$ be nilpotent. We investigate the
centralizer $\mathfrak{g}^{e}=\{x\in\mathfrak{g}:[e,x]=0\}$ of $e$
in $\mathfrak{g}$ and the centre of centralizer $\mathfrak{z}(\mathfrak{g}^{e})=\{x\in\mathfrak{g}^{e}:[x,y]=0\text{ for all }y\in\mathfrak{g}^{e}\}$
of $e$ in $\mathfrak{g}$. A lot of research has been done on the
centralizer and the centre of centralizer of nilpotent elements in
the theory of Lie algebras. Although there are similarities between
the theory of Lie superalgebras and the theory of Lie algebras, there
is a lot less study in this direction in the case of Lie superalgebras
and the structural theory of nilpotent orbits in Lie superalgebras
remains to be better understood. In this paper, we calculate bases
for $\mathfrak{g}^{e}$ and $\mathfrak{z}(\mathfrak{g}^{e})$ and
study various properties relating $e$ with $\mathfrak{g}^{e}$ and
$\mathfrak{z}(\mathfrak{g}^{e})$. 

Research on the centralizer of nilpotent elements and their centres
in the case of Lie algebras has been intensively developed since Springer
\cite{Springer1966} considered the centralizer $G^{u}$ of a unipotent
element $u$ in a simple algebraic group $G$. Many mathematicians
undertook further study of $G^{u}$, the reader is referred to the
introduction of \cite{Lawther2008} for an overview of research on
$G^{u}$. For classical Lie algebras over an algebraically closed
field of arbitrary characteristic, Jantzen gave an explicit account
of the structure of $\mathfrak{g}^{e}$ in \cite{Janzten} and Yakimova
worked out bases for $\mathfrak{z}(\mathfrak{g}^{e})$ in \cite{Yakimova2009}.
In \cite{Lawther2008}, Lawther--Testerman dealt with the centralizer
$G^{u}$ and its centre $Z(G^{u})$ over a field of characteristic
$0$ or a good prime based on Yakimova's results. In \cite{han-classical,han-exp},
the author identified $\mathfrak{g}^{e}$ and $\mathfrak{z}(\mathfrak{g}^{e})$
for basic classical Lie superalgebras over a field of characteristic
zero and obtained analog of results of Lawther--Testerman \cite{Lawther2008}
for those Lie superalgebras.

Define $\mathfrak{g}_{\mathbb{C}}$ to be a finite-dimensional basic
classical Lie superalgebra over $\mathbb{C}$ and write $\Phi$ for
a root system of $\mathfrak{g}_{\mathbb{C}}$. By \cite[Theorem 3.9]{Iohara=000026Koga},
there exists a Chevalley basis $\mathfrak{B}=\{e_{\alpha}:\alpha\in\Phi\}\cup\{h_{i}:1\leq i\leq s\}$
of $\mathfrak{g}_{\mathbb{C}}$ such that $[h_{i},e_{\alpha}]=\langle\alpha,\alpha_{i}\rangle e_{\alpha}$
and $[e_{\alpha},e_{\beta}]=N_{\alpha,\beta}e_{\alpha+\beta}$ where
$\langle\alpha,\alpha_{i}\rangle$ is defined in (\ref{eq:<beta,alpha>})
and $N_{\alpha,\beta}\in\mathbb{Z}$ can be determined explicitly.
Let $\mathfrak{g}_{\mathbb{Z}}\subseteq\mathfrak{g}_{\mathbb{C}}$
be the Chevalley $\mathbb{Z}$-form of $\mathfrak{g}_{\mathbb{C}}$,
i.e. $\mathfrak{g}_{\mathbb{Z}}$ is the $\mathbb{Z}$-span of $\mathfrak{B}$.
Denote by $\mathfrak{h}_{\mathbb{Z}}=\langle h_{i}:1\leq i\leq s\rangle_{\mathbb{Z}}$
a Cartan subalgebra of $\mathfrak{g}_{\mathbb{Z}}$. We can view $\mathfrak{g}=\mathfrak{g}_{\mathbb{K}}$
where $\mathfrak{g}_{\mathbb{K}}=\mathbb{K}\otimes_{\mathbb{Z}}\mathfrak{g}_{\mathbb{Z}}$.
It is natural to ask what is the stucture of $\mathfrak{g}^{e}$ and
$\mathfrak{z}(\mathfrak{g}^{e})$ in case $p>0$ is a good prime for
$\mathfrak{g}$. For $\mathfrak{g}=\mathfrak{gl}(m|n)$, the construction
of $\mathfrak{g}^{e}$ over a field of prime characteristic and that
of zero characteristic are identical by \cite{WANGZHAO} and \cite{Hoyt2012}.
However, the structure of $\mathfrak{g}^{e}$ in case of $\mathfrak{g}=\mathfrak{sl}(m|n)$
for $m\neq n$, $\mathfrak{psl}(n|n)$, $\mathfrak{osp}(m|2n)$ and
three exceptional types have not been considered yet. In this paper,
we aim to give a description of $\mathfrak{g}^{e}$ and further calculate
the dimension of $\mathfrak{z}(\mathfrak{g}^{e})$ for the above types
of Lie superalgebras. 

In the remaining part of this introduction, we give a more detailed
survey of our results.

Fix $\mathbb{K}$ an algebraically closed field with a good prime
characteristic $p>0$, see Definition \ref{def:good prime}. Let $\mathfrak{g}=\mathfrak{g}_{\mathbb{K}}=\mathfrak{g}_{\bar{0}}\oplus\mathfrak{g}_{\bar{1}}$
be one of the Lie superalgebras in Table \ref{tab:Basic-classical}.
Let $G_{\bar{0}}$ be the reductive algebraic group over $\mathbb{K}$
given as in Table \ref{tab:algebraic-groups} such that $\mathrm{Lie}(G_{\bar{0}})=\mathfrak{g}_{\bar{0}}$.
Then there is a representation $\rho:G_{\bar{0}}\rightarrow\mathrm{GL}(\mathfrak{g}_{\bar{1}})$
such that $d_{\rho}:\mathrm{Lie}(G_{\bar{0}})\rightarrow\mathfrak{gl}(\mathfrak{g}_{\bar{1}})$
determines the adjoint action of $\mathfrak{g}_{\bar{0}}$ on $\mathfrak{g}_{\bar{1}}$.

\begin{doublespace}
\noindent %
\begin{longtable}[c]{|>{\centering}p{3cm}|c|}
\caption{\label{tab:algebraic-groups}Algebraic groups $G_{\bar{0}}$}
\tabularnewline
\endfirsthead
\hline 
Lie superalgebras $\mathfrak{g}$ & Algebraic groups $G_{\bar{0}}$\tabularnewline
\hline 
\hline 
$\mathfrak{sl}(m|n),m\neq n$ & $\left\{ (A,B)\in\mathrm{GL}_{m}(\mathbb{K})\times\mathrm{GL}_{n}(\mathbb{K}):\mathrm{det}(A)=\mathrm{det}(B)\right\} $\tabularnewline
\hline 
$\mathfrak{psl}(n|n)$ & $\left\{ (A,B)\in\mathrm{GL}_{n}(\mathbb{K})\times\mathrm{GL}_{n}(\mathbb{K}):\mathrm{det}(A)=\mathrm{det}(B)\}/\left\{ aI_{n|n}:a\in\mathbb{K}^{\times}\right\} \right\} $\tabularnewline
\hline 
$\mathfrak{osp}(m|2n)$ & $\mathrm{O}_{m}(\mathbb{K})\times\mathrm{Sp}_{2n}(\mathbb{K})$\tabularnewline
\hline 
$D(2,1;\alpha)$ & $\mathrm{SL}_{2}(\mathbb{K})\times\mathrm{SL}_{2}(\mathbb{K})\times\mathrm{SL}_{2}(\mathbb{K})$\tabularnewline
\hline 
$G(3)$ & $\mathrm{SL}_{2}(\mathbb{K})\times G_{2}$\tabularnewline
\hline 
$F(4)$ & $\mathrm{SL}_{2}(\mathbb{K})\times\mathrm{Spin}_{7}(\mathbb{K})$\tabularnewline
\hline 
\end{longtable}
\end{doublespace}

For each nilpotent element $e\in(\mathfrak{g}_{\mathbb{C}})_{\bar{0}}$,
the Jacobson--Morozov Theorem allows one to associate an $\mathfrak{sl}_{2}$-triple
$\{e,h,f\}\subseteq(\mathfrak{g}_{\mathbb{C}})_{\bar{0}}$ to $e$.
According to \cite[Section 3]{Iohara=000026Koga}, the $\mathfrak{sl}_{2}$-triple
can be chosen such that $\{e,h,f\}\subseteq(\mathfrak{g}_{\mathbb{Z}})_{\bar{0}}$
where $h\in\mathfrak{h}_{\mathbb{Z}}$ is of the form $h=\sum_{i=1}^{s}c_{i}h_{i}$
for $c_{i}\in\mathbb{Z}$ and $e=\sum_{\alpha\in\Phi}e_{\alpha}$.
Note that the $\mathrm{ad}h$-grading of $\mathfrak{g}_{\mathbb{\mathbb{Z}}}=\bigoplus_{j\in\mathbb{Z}}\mathfrak{g}_{\mathbb{Z}}(j;\mathrm{ad}h)$
is given by $\mathfrak{g}_{\mathbb{Z}}(j;\mathrm{ad}h)=\{x\in\mathfrak{g}_{\mathbb{Z}}:[h,x]=jx\}$
and this grading can be extended to $\mathfrak{g}_{\mathbb{C}}$.
By the representation theory of $\mathfrak{sl}_{2}(\mathbb{C})$,
we can determine $\mathrm{ad}h$-eigenvalues of elements of $\mathfrak{g}_{\mathbb{C}}^{e}$.
Based on the choice of $e$, we also can view $e\in\mathfrak{g}_{\mathbb{K}}$.
We calculate bases for $\mathfrak{g}^{e}$ and $\mathfrak{z}(\mathfrak{g}^{e})$
for Lie superalgebras of type $A(m,n)$, $B(m,n)$, $C(n)$, $D(m,n)$
and $D(2,1;\alpha)$, $G(3)$, $F(4)$ in Sections \ref{sec:A(m,n)}--\ref{sec:osp}
and \ref{sec:exceptional} respectively. In particular, we have the
following result. 
\begin{thm}
\label{thm:Main_Thm_1}There exists a basis $\mathfrak{B}^{e}\subseteq\mathfrak{g}_{\mathbb{Z}}$
of $\mathfrak{g}_{\mathbb{C}}^{e}$ such that $\mathfrak{B}^{e}$
viewed in $\mathfrak{g}^{e}$ is a basis for $\mathfrak{g}^{e}$.
Similarly we can find a basis $\mathfrak{B}_{\mathfrak{z}}^{e}\subseteq\mathfrak{g}_{\mathbb{Z}}$
of $\mathfrak{z}(\mathfrak{g}_{\mathbb{C}}^{e})$ such that $\mathfrak{B}_{\mathfrak{z}}^{e}$
viewed in $\mathfrak{z}(\mathfrak{g}^{e})$ is a basis for $\mathfrak{z}(\mathfrak{g}^{e})$.
Note that when $\mathfrak{g}=\mathfrak{sl}(m|n)$ for $m\neq n$ or
$\mathfrak{psl}(n|n)$ for $n>1$, we require that $\mathrm{char}(\mathbb{K})=p$
does not divide $m$ and $n$.
\end{thm}

Note that in good characteristic there is a substitute for $\mathfrak{sl}{}_{2}$-triples,
so called associated cocharacters, see Definition \ref{def:cocharacter}
below. Let $\tau:\mathbb{K}^{\times}\rightarrow G_{\bar{0}}$ be a
cocharacter associated to $e$. Denote by $\mathfrak{g}=\bigoplus_{j\in\mathbb{Z}}\mathfrak{g}(j;\tau)$
the $\tau$-grading on $\mathfrak{g}$ where $\mathfrak{g}(j;\tau)=\{x\in\mathfrak{g}:\mathrm{Ad}(\tau(t))(x)=t^{j}x\text{ for all }t\in\mathbb{K}^{\times}\}$.
In \cite[Section 3]{WANGZHAO}, Wang--Zhao studied properties of
the $\tau$-grading on $\mathfrak{g}$ with some restrictions on $p$.
Combining Theorem \ref{thm:Main_Thm_1} with Lemma \ref{lem:g(j)},
we obtain the following theorem which gives a more general statement
on the $\tau$-grading on $\mathfrak{g}$.
\begin{thm}
\label{thm:Main_Thm_2}Let $\mathfrak{g}=\mathfrak{g}_{\bar{0}}\oplus\mathfrak{g}_{\bar{1}}$
be one of the basic classical Lie superalgebras in Table \ref{tab:Basic-classical}
and $e\in\mathfrak{g}_{\bar{0}}$ be nilpotent. Then the cocharacter
$\tau:\mathbb{K}^{\times}\rightarrow G_{\bar{0}}$ associated to $e$
defines a $\mathbb{Z}$-grading $\mathfrak{g}=\bigoplus_{j\in\mathbb{Z}}\mathfrak{g}(j;\tau)$
such that $\mathfrak{g}^{e}=\bigoplus_{j\geq0}\mathfrak{g}^{e}(j;\tau)$
and $\dim\mathfrak{g}^{e}(j;\tau)=\dim\mathfrak{g}(j;\tau)-\dim\mathfrak{g}((j+2);\tau)$
for $j\geq0$. 
\end{thm}

Our next result focuses on the $\mathbb{Z}$-grading on $\mathfrak{z}(\mathfrak{g}^{e})$.
In the case of simple Lie algebras, the nilpotent element $e$ spans
the degree $2$ part of the $\mathrm{ad}h$-grading of the centre
of centralizer of $e$. It is known as the Brylinsky--Kostant theorem,
see \cite{Brylinski1989,Kostant1963}. Write $\mathfrak{z}(j;\mathrm{ad}h)=\mathfrak{z}(\mathfrak{g}_{\mathbb{C}}^{e})\cap\mathfrak{g}(j;\mathrm{ad}h)$
and $\mathfrak{z}(j;\tau)=\mathfrak{z}(\mathfrak{g}^{e})\cap\mathfrak{g}(j;\tau)$.
Theorem \ref{thm:Main_Thm_3} can be viewed as the Lie superalgebra
version of the Brylinsky--Kostant theorem.
\begin{thm}
\label{thm:Main_Thm_3}Let $\mathfrak{g}=\mathfrak{g}_{\bar{0}}\oplus\mathfrak{g}_{\bar{1}}$
be one of the basic classical Lie superalgebras in Table \ref{tab:Basic-classical}
and $e\in\mathfrak{g}_{\bar{0}}$ be nilpotent. Then $\mathfrak{z}(2;\mathrm{ad}h)=\langle e\rangle$
and $\mathfrak{z}(2;\tau)=\langle e\rangle$.
\end{thm}

The element $e$ is called \textit{reachable} if $e\in[\mathfrak{g}^{e},\mathfrak{g}^{e}]$,
such element were first considered by Elashvili and Grélaud and called
\textit{compact} in \cite{Elashvili=000026Gr=0000E9laud}. The element
$e$ is called \textit{strongly reachable} if $\mathfrak{g}^{e}=[\mathfrak{g}^{e},\mathfrak{g}^{e}]$.
We say that $e$ satisfies the \textit{Panyushev property} \cite{Panyushev}
if in the $\tau$-grading $\mathfrak{g}^{e}=\bigoplus_{j\geq0}\mathfrak{g}(j)$,
the subalgebra $\mathfrak{g}^{e}(\geq1)=\bigoplus_{j\geq1}\mathfrak{g}(j)$
is generated by $\mathfrak{g}^{e}(1)$. In the case of Lie algebras,
Panyushev \cite{Panyushev} showed that for $\mathfrak{g}$ of type
$A_{n}$, a nilpotent element $e$ is reachable if and only if $e$
satisfies the Panyushev property. The result was extended to cases
for $\mathfrak{g}$ of type $B_{n}$, $C_{n}$, $D_{n}$ by Yakimova
\cite{Yakimova2010} and for $\mathfrak{g}$ of exceptional types
by de Graaf \cite{de Graaf2013}. The author \cite{han-reachable}
gave the classification of even elements that are reachable, strongly
reachable or satisfying the Panyushev property in basic classical
Lie superalgebras over $\mathbb{C}$. Our final result extends the
results in \cite{han-reachable} to good characteristic, which illustrates
the relation between the property of being reachable and the Panyushev
property.
\begin{thm}
Let $\mathfrak{g}=\mathfrak{g}_{\bar{0}}\oplus\mathfrak{g}_{\bar{1}}$
be one of exceptional Lie superalgebras $D(2,1;\alpha)$, $G(3)$
, $F(4)$ and $e\in\mathfrak{g}_{\bar{0}}$ be nilpotent. Then $e$
is reachable if and only if $e$ satisfies the Panyushev property
except for $\mathfrak{g}=G(3)$ and $e=x_{1}$ or $e=E+x_{1}$. The
nilpotent orbits that are reachable, strongly reachable or satisfying
the Panyushev property are listed in Tables \ref{tab:Reachable D}--\ref{tab:Reachable F}.
\end{thm}

This paper is organized as follows. We first recall some fundamental
concepts of Lie superalgebras such as basic classical Lie superalgebras,
root systems and cocharacters associated to nilpotent elements in
Section \ref{sec:Preliminaries}. For each system of positive roots,
we identify the highest root in it. In Sections \ref{sec:A(m,n)}--\ref{sec:osp},
we determine bases of $\mathfrak{g}^{e}$ and $\mathfrak{z}(\mathfrak{g}^{e})$
for $\mathfrak{g}=\mathfrak{sl}(m|n)$ for $m\neq n$, $\mathfrak{psl}(n|n)$
and $\mathfrak{osp}(m|2n)$. Bases of $\mathfrak{g}^{e}$ and $\mathfrak{z}(\mathfrak{g}^{e})$
for an exceptional Lie superalgebra $\mathfrak{g}$ are given in Section
\ref{sec:exceptional}. Using the structure of $\mathfrak{g}^{e}$
for an exceptional Lie superalgebra $\mathfrak{g}$ in Section \ref{sec:exceptional},
we determine which even nilpotent elements are reachable, strongly
reachable or satisfy the Panyushev property in Section \ref{sec:Reachability}. 

\noindent \textbf{Acknowledgements.} The author acknowledges financial
support from the Engineering and Physical Sciences Research Council
(EP/W522478/1). We would like to thank Simon Goodwin for very useful
discussions on the subject of this paper.

\section{Preliminaries and notations\label{sec:Preliminaries}}

\subsection{Basic classical Lie superalgebras}

Finite-dimensional simple Lie superalgebras over $\mathbb{C}$ were
classified by Kac in \cite{Kac}. Among those simple Lie superalgebras,
we focus on basic classical Lie superalgebras in this paper. Recall
that a finite-dimensional simple Lie superalgebra $\mathfrak{g}=\mathfrak{g}_{\bar{0}}\oplus\mathfrak{g}_{\bar{1}}$
is called \textit{a basic classical Lie superalgebra} if the even
part $\mathfrak{g}_{\bar{0}}$ is a reductive Lie algebra and there
exists a non-degenerate supersymmetric invariant even bilinear form
$(\cdotp,\cdotp)$ on $\mathfrak{g}$. Note that these Lie superalgebras
are well defined over $\mathbb{K}$ and remain to simple by \cite[Section 2]{WANGZHAO}.
Below in Table \ref{tab:Basic-classical} we recall the list of basic
classical Lie superalgebras over $\mathbb{K}$ that are not Lie algebras,
they are $A(m,n)$ for $m\neq n$, $A(n,n)$ for $n>1$, $B(m,n)$,
$C(n)$, $D(m,n)$ and three exceptional types $D(2,1;\alpha)$, $G(3)$,
$F(4)$. The explicit construction of $A(m,n)$, $A(n,n)$, $B(m,n)$,
$C(n)$, $D(m,n)$ can be found for example in \cite[Sections 3--4]{han-classical}
and that of $D(2,1;\alpha)$, $G(3)$, $F(4)$ can be found for example
in \cite[Sections 4--6]{han-exp}.

\begin{table}[H]
\begin{centering}
\begin{tabular}{|c|c|}
\hline 
$\mathfrak{g}=\mathfrak{g}_{\bar{0}}\oplus\mathfrak{g}_{\bar{1}}$ & $\mathfrak{g}_{\bar{0}}$\tabularnewline
\hline 
\hline 
$A(m,n)=\mathfrak{sl}(m|n)$, $m,n\geq1$, $m\neq n$ & $\mathfrak{sl}_{m}(\mathbb{K})\oplus\mathfrak{sl}_{n}(\mathbb{K})\oplus\mathbb{K}$\tabularnewline
\hline 
$A(n,n)=\mathfrak{psl}(n|n)$, $n>1$ & $\mathfrak{sl}_{n}(\mathbb{K})\oplus\mathfrak{sl}_{n}(\mathbb{K})$\tabularnewline
\hline 
$B(m,n)=\mathfrak{osp}(m|2n)$, $m$ is odd, $m,n\geq1$ & $\mathfrak{o}_{m}(\mathbb{K})\oplus\mathfrak{sp}_{2n}(\mathbb{K})$\tabularnewline
\hline 
$C(n)=\mathfrak{osp}(2|2n)$, $n\geq1$ & $\mathfrak{o}_{2}(\mathbb{K})\oplus\mathfrak{sp}_{2n}(\mathbb{K})$\tabularnewline
\hline 
$D(m,n)=\mathfrak{osp}(m|2n)$, $m\geq4$ is even, $n\geq1$ & $\mathfrak{o}_{m}(\mathbb{K})\oplus\mathfrak{sp}_{2n}(\mathbb{K})$\tabularnewline
\hline 
$D(2,1;\alpha)$, $\alpha\in\mathbb{K}\backslash\{0,-1\}$ & $\mathfrak{sl}_{2}(\mathbb{K})\oplus\mathfrak{sl}_{2}(\mathbb{K})\oplus\mathfrak{sl}_{2}(\mathbb{K})$\tabularnewline
\hline 
$G(3)$ & $\mathfrak{sl}_{2}(\mathbb{K})\oplus G_{2}$\tabularnewline
\hline 
$F(4)$ & $\mathfrak{sl}_{2}(\mathbb{K})\oplus\mathfrak{so}_{7}(\mathbb{K})$\tabularnewline
\hline 
\end{tabular}
\par\end{centering}
\caption{\label{tab:Basic-classical}Basic classical Lie superalgebras}
\end{table}

\subsection{Systems of positive roots and the highest root\label{sec:highest root}}

\noindent Let $\mathfrak{g}$ be one of the basic classical Lie superalgebras
in Table \ref{tab:Basic-classical} and let $\Phi=\Phi_{\bar{0}}\cup\Phi_{\bar{1}}$
be a root system for $\mathfrak{g}$. Let $\Phi^{+}$ be a system
of positive roots and $\Pi=\{\alpha_{1},\dots,\alpha_{l}\}$ be the
corresponding system of simple roots. Based on results in \cite{han-classical,han-exp},
there is in general more than one system of simple roots. For each
system of positive roots, there is an unique \textit{highest root}
$\widetilde{\alpha}$ in $\Phi^{+}$ such that $\widetilde{\alpha}=\sum_{i=1}^{l}a_{i}\alpha_{i}$
and for any other $\beta=\sum_{i=1}^{l}b_{i}\alpha_{i}\in\Phi^{+}$,
we have $b_{i}\leq a_{i}$ for all $i$. 

In the following part of this subsection, we adopt notations in \cite[Sections 3--4]{han-classical}
and \cite[Sections 4--6]{han-exp} for root systems of $\mathfrak{g}$
and systems of simple roots. In particular, we use the results in
\cite{han-classical,han-exp} to determine the highest root in $\Phi^{+}$
for each system of simple roots.

\noindent \textbf{1. $A(m,n)$}

Let $V=V_{\bar{0}}\oplus V_{\bar{1}}$ be a finite-dimensional vector
space over $\mathbb{K}$ such that $\dim V_{\bar{0}}=m$ and $\dim V_{\bar{1}}=n$.
For a homogeneous basis $\{v_{1},\dots,v_{m+n}\}$ of $V$, we denote
by $\eta_{i}$ the parity of $v_{i}$, i.e. $\eta_{i}=\bar{0}$ if
$v_{i}\in V_{\bar{0}}$ and $\eta_{i}=\bar{1}$ if $v_{i}\in V_{\bar{1}}$.
Let $e_{ij}$ be the $ij$-matrix unit. A Cartan subalgebra $\mathfrak{h}$
of $\mathfrak{gl}(m|n)$ has a basis $\{h_{i}=e_{i,i}:1\leq i\leq m+n\}$
and a dual basis $\{\varepsilon_{i}^{\eta_{i}}\in\mathfrak{h}^{*}:1\leq i\leq m+n\}$
is defined by $\varepsilon_{i}^{\eta_{i}}(h_{j})=\delta_{ij}$ such
that $(\varepsilon_{i}^{\eta_{i}},\varepsilon_{j}^{\eta_{j}})=(-1)^{\eta_{i}}\delta_{ij}$.
Then $\mathfrak{g}=\mathfrak{sl}(V)$ or $\mathfrak{psl}(V)$ has
a root system $\Phi=\Phi_{\bar{0}}\cup\Phi_{\bar{1}}$ such that $\Phi_{\bar{0}}=\{\varepsilon_{i}^{\eta_{i}}-\varepsilon_{j}^{\eta_{j}}:1\leq i,j\leq m+n,i\neq j,\eta_{i}=\eta_{j}\}$
and $\Phi_{\bar{1}}=\{\varepsilon_{i}^{\eta_{i}}-\varepsilon_{j}^{\eta_{j}}:1\leq i,j\leq m+n,i\neq j,\eta_{i}\neq\eta_{j}\}$.
The system of simple roots is given by $\Pi=\{\alpha_{i}=\varepsilon_{i}^{\eta_{i}}-\varepsilon_{i+1}^{\eta_{i+1}}:1\leq i<m+n\}$
and the highest root is $\widetilde{\alpha}=\varepsilon_{1}^{\eta_{1}}-\varepsilon_{m+n}^{\eta_{m+n}}=\sum_{i=1}^{m+n-1}\alpha_{i}$.

\noindent \textbf{2. $B(m,n)$, $C(n)$ or $D(m,n)$} 

Let $V=V_{\bar{0}}\oplus V_{\bar{1}}$ be a finite-dimensional vector
space over $\mathbb{K}$ such that $\dim V_{\bar{0}}=m$ and $\dim V_{\bar{1}}=2n$
and let $\mathfrak{g}=\mathfrak{osp}(V)$. Write $l=\left\lfloor \frac{m}{2}\right\rfloor $.
A Cartan subalgebra $\mathfrak{h}$ of $\mathfrak{g}$ has a basis
$\{h_{i}=e_{i,i}-e_{-i,-i}:1\leq i\leq l+n\}$ and a dual basis $\{\varepsilon_{i}^{\eta_{i}}\in\mathfrak{h}^{*}:1\leq i\leq l+n\}$
is defined by $\varepsilon_{i}^{\eta_{i}}(h_{j})=\delta_{ij}$ such
that $(\varepsilon_{i}^{\eta_{i}},\varepsilon_{j}^{\eta_{j}})=(-1)^{\eta_{i}}\delta_{ij}$. 

When $m$ is odd, a root system for $\mathfrak{g}$ is $\Phi_{\bar{0}}=\{\pm\varepsilon_{i}^{\eta_{i}}\pm\varepsilon_{j}^{\eta_{j}}:i\neq j,\eta_{i}=\eta_{j}\}\cup\{\pm\varepsilon_{i}^{\bar{0}}\}\cup\{\pm2\varepsilon_{i}^{\bar{1}}\}$
and $\Phi_{\bar{1}}=\{\pm\varepsilon_{i}^{\bar{1}}\}\cup\{\pm\varepsilon_{i}^{\eta_{i}}\pm\varepsilon_{j}^{\eta_{j}}:i\neq j,\eta_{i}\neq\eta_{j}\}$.
A system of simple roots is $\{\alpha_{i}=\varepsilon_{i}^{\eta_{i}}-\varepsilon_{i+1}^{\eta_{i+1}}:1\leq i\leq l+n-1\}\cup\{\alpha_{l+n}=\varepsilon_{l+n}^{\eta_{l+n}}\}$.
The associated system of positive roots is $\Phi^{+}=\{\varepsilon_{i}^{\eta_{i}}\pm\varepsilon_{j}^{\eta_{j}},\varepsilon_{k}^{\eta_{k}},2\varepsilon_{t}^{\bar{1}}\}$,
where the highest root is 
\[
\widetilde{\alpha}=\begin{cases}
\varepsilon_{1}^{\eta_{1}}+\varepsilon_{2}^{\eta_{2}}=\alpha_{1}+2\sum_{i=2}^{l+n}\alpha_{i} & \text{if }\eta_{1}=\bar{0};\\
2\varepsilon_{1}^{\bar{1}}=2\sum_{i=1}^{l+n}\alpha_{i} & \text{if }\eta_{1}=\bar{1}.
\end{cases}
\]

When $m$ is even, a root system for $\mathfrak{g}$ is $\Phi_{\bar{0}}=\{\pm\varepsilon_{i}^{\eta_{i}}\pm\varepsilon_{j}^{\eta_{j}}:i\neq j,\eta_{i}=\eta_{j}\}\cup\{\pm2\varepsilon_{i}^{\bar{1}}\}$
and $\Phi_{\bar{1}}=\{\pm\varepsilon_{i}^{\eta_{i}}\pm\varepsilon_{j}^{\eta_{j}}:i\neq j,\eta_{i}\neq\eta_{j}\}$.
Note that there are two possibilities for the systems of simple roots
with the same associated system of positive roots $\Phi^{+}=\{\varepsilon_{i}^{\eta_{i}}\pm\varepsilon_{j}^{\eta_{j}},2\varepsilon_{t}^{\bar{1}}\}$.
There are 

$\bullet$ $\Pi_{1}=\{\alpha_{i}=\varepsilon_{i}^{\eta_{i}}-\varepsilon_{i+1}^{\eta_{i+1}}:1\leq i\leq l+n-2\}\cup\{\alpha_{l+n-1}=\varepsilon_{l+n-1}^{\eta_{l+n-1}}-\varepsilon_{l+n}^{\bar{1}}\}\cup\{\alpha_{l+n}=2\varepsilon_{l+n}^{\bar{1}}\}$,
and the highest root is 
\[
\widetilde{\alpha}_{(1)}=\begin{cases}
\varepsilon_{1}^{\eta_{1}}+\varepsilon_{2}^{\eta_{2}}=\alpha_{1}+2\sum_{i=2}^{l+n-1}\alpha_{i}+\alpha_{l+n} & \text{if }\eta_{1}=\bar{0};\\
2\varepsilon_{1}^{\bar{1}}=2\sum_{i=1}^{l+n-1}\alpha_{i}+\alpha_{l+n} & \text{if }\eta_{1}=\bar{1}.
\end{cases}
\]

$\bullet$ $\Pi_{2}=\{\alpha_{i}=\varepsilon_{i}^{\eta_{i}}-\varepsilon_{i+1}^{\eta_{i+1}}:1\leq i\leq l+n-2\}\cup\{\alpha_{l+n-1}=\varepsilon_{l+n-1}^{\eta_{l+n-1}}-\varepsilon_{l+n}^{\eta_{l+n}}\}\cup\{\alpha_{l+n}=\varepsilon_{l+n-1}^{\eta_{l+n-1}}+\varepsilon_{l+n}^{\eta_{l+n}}\}$,
and the highest root is 
\[
\widetilde{\alpha}_{(2)}=\begin{cases}
\varepsilon_{1}^{\eta_{1}}+\varepsilon_{2}^{\eta_{2}}=\alpha_{1}+2\sum_{i=2}^{l+n-2}\alpha_{i}+\alpha_{l+n-1}+\alpha_{l+n} & \text{if }\eta_{1}=\bar{0};\\
2\varepsilon_{1}^{\bar{1}}=2\sum_{i=1}^{l+n-2}\alpha_{i}+\alpha_{l+n-1}+\alpha_{l+n} & \text{if }\eta_{1}=\bar{1}.
\end{cases}
\]

\noindent \textbf{3. $D(2,1;\alpha)$ (with $\alpha\neq0,-1$) }

The Lie superalgebra $\mathfrak{g}=D(2,1;\alpha)$ has a root system
$\Phi=\Phi_{\bar{0}}\cup\Phi_{\bar{1}}$ such that $\Phi_{\bar{0}}=\lbrace\pm2\beta_{1},\pm2\beta_{2},\pm2\beta_{3}\rbrace$,
$\Phi_{\bar{1}}=\lbrace\pm\beta_{1}\pm\beta_{2}\pm\beta_{3}\rbrace$
where $\{\beta_{1},\beta_{2},\beta_{3}\}$ is an orthogonal basis
of $\mathbb{R}\Phi$ such that $(\beta_{1},\beta_{1})=\frac{1}{2}$,
$(\beta_{2},\beta_{2})=-\frac{1}{2}\alpha-\frac{1}{2}$ and $(\beta_{3},\beta_{3})=\frac{1}{2}\alpha$.
According to \cite[Subsection 4.2]{han-exp}, there are four conjugacy
classes of systems of simple roots which we list in the table below.

\begin{doublespace}
\noindent %
\begin{longtable}[c]{|>{\centering}m{4cm}|>{\centering}m{5cm}|>{\centering}m{3.5cm}|}
\caption{\label{tab:highest root D}The highest root in system of positive
roots for $D(2,1;\alpha)$}
\tabularnewline
\endfirsthead
\hline 
System of simple roots $\Pi=\{\alpha_{1},\alpha_{2},\alpha_{3}\}$ & Associated system of positive roots $\Phi^{+}$ & The highest root $\widetilde{\alpha}$\tabularnewline
\hline 
\hline 
$\{2\beta_{1},-\beta_{1}+\beta_{2}-\beta_{3},2\beta_{3}\rbrace$ & $\lbrace2\beta_{i},\pm\beta_{1}+\beta_{2}\pm\beta_{3}:i=1,2,3\rbrace$ & $2\beta_{2}=\alpha_{1}+2\alpha_{2}+\alpha_{3}$\tabularnewline
\hline 
$\{2\beta_{1},-\beta_{1}-\beta_{2}+\beta_{3},2\beta_{2}\rbrace$ & $\lbrace2\beta_{i},\pm\beta_{1}\pm\beta_{2}+\beta_{3}:i=1,2,3\rbrace$ & $2\beta_{3}=\alpha_{1}+2\alpha_{2}+\alpha_{3}$\tabularnewline
\hline 
$\{2\beta_{3},\beta_{1}-\beta_{2}-\beta_{3},2\beta_{2}\rbrace$ & $\lbrace2\beta_{i},\beta_{1}\pm\beta_{2}\pm\beta_{3}:i=1,2,3\rbrace$ & $2\beta_{1}=\alpha_{1}+2\alpha_{2}+\alpha_{3}$\tabularnewline
\hline 
$\{-\beta_{1}+\beta_{2}+\beta_{3},\beta_{1}-\beta_{2}+\beta_{3},\beta_{1}+\beta_{2}-\beta_{3}\rbrace$ & $\lbrace2\beta_{i},\pm\beta_{1}+\beta_{2}+\beta_{3},\beta_{1}-\beta_{2}+\beta_{3},\beta_{1}+\beta_{2}-\beta_{3}:i=1,2,3\rbrace$ & $\beta_{1}+\beta_{2}+\beta_{3}=\alpha_{1}+\alpha_{2}+\alpha_{3}$\tabularnewline
\hline 
\end{longtable}
\end{doublespace}

\noindent \textbf{4. $G(3)$ }

The Lie superalgebra $\mathfrak{g}=G(3)$ has the root system $\Phi=\Phi_{\bar{0}}\cup\Phi_{\bar{1}}$
such that $\Phi_{\bar{0}}=\{\pm2\delta,\varepsilon_{i}-\varepsilon_{j},\pm\varepsilon_{i}:1\leq i\neq j\leq3\}$
and $\Phi_{\bar{1}}=\{\pm\delta\pm\varepsilon_{i},\pm\delta:1\leq i\leq3\}$
where $\{\delta,\varepsilon_{1,}\varepsilon_{2},\varepsilon_{3}\}$
are elements of $(\bigoplus_{i}\mathbb{C}\varepsilon_{i}\oplus\mathbb{C}\delta)/\mathbb{C}(\varepsilon_{1}+\varepsilon_{2}+\varepsilon_{3})$
such that $(\delta,\delta)=2$, $(\varepsilon_{i},\varepsilon_{j})=1-3\delta_{ij}$,
and $(\delta,\varepsilon_{i})=0$. According to \cite[Subsection 5.2]{han-exp},
there are four conjugacy classes of systems of simple roots which
we list in the table below.

\begin{doublespace}
\noindent %
\begin{longtable}[c]{|>{\centering}m{3.5cm}|>{\centering}m{5cm}|>{\centering}m{4cm}|}
\caption{\label{tab:highest root G}The highest root in system of positive
roots for $G(3)$}
\tabularnewline
\endfirsthead
\hline 
System of simple roots $\Pi=\{\alpha_{1},\alpha_{2},\alpha_{3}\}$ & Associated system of positive roots $\Phi^{+}$ & The highest root $\widetilde{\alpha}$\tabularnewline
\hline 
\hline 
$\{\delta+\varepsilon_{3},\varepsilon_{1},\varepsilon_{2}-\varepsilon_{1}\rbrace$ & $\lbrace\varepsilon_{1},\varepsilon_{2},-\varepsilon_{3},\varepsilon_{2}-\varepsilon_{1},\varepsilon_{1}-\varepsilon_{3},\varepsilon_{2}-\varepsilon_{3},\delta,2\delta,\delta\pm\varepsilon_{i}:i=1,2,3\rbrace$ & $2\delta=2\alpha_{1}+4\alpha_{2}+2\alpha_{3}$\tabularnewline
\hline 
$\{-\delta-\varepsilon_{3},\delta-\varepsilon_{2},\varepsilon_{2}-\varepsilon_{1}\rbrace$ & $\lbrace\varepsilon_{1},\varepsilon_{2},-\varepsilon_{3},\varepsilon_{2}-\varepsilon_{1},\varepsilon_{1}-\varepsilon_{3},\varepsilon_{2}-\varepsilon_{3},\delta,2\delta,\delta\pm\varepsilon_{1},\delta\pm\varepsilon_{2},\pm\delta-\varepsilon_{3}\rbrace$ & $\delta-\varepsilon_{3}=3\alpha_{1}+4\alpha_{2}+2\alpha_{3}$\tabularnewline
\hline 
$\{\delta,-\delta+\varepsilon_{1},\varepsilon_{2}-\varepsilon_{1}\rbrace$ & $\lbrace\varepsilon_{1},\varepsilon_{2},-\varepsilon_{3},\varepsilon_{2}-\varepsilon_{1},\varepsilon_{1}-\varepsilon_{3},\varepsilon_{2}-\varepsilon_{3},\delta,2\delta,\pm\delta+\varepsilon_{1},\pm\delta+\varepsilon_{2},\pm\delta-\varepsilon_{3}\rbrace$ & $\varepsilon_{2}-\varepsilon_{3}=3\alpha_{1}+2\alpha_{2}+2\alpha_{3}$\tabularnewline
\hline 
$\{\varepsilon_{1},-\delta+\varepsilon_{2},\delta-\varepsilon_{1}\rbrace$ & $\lbrace\varepsilon_{1},\varepsilon_{2},-\varepsilon_{3},\varepsilon_{2}-\varepsilon_{1},\varepsilon_{1}-\varepsilon_{3},\varepsilon_{2}-\varepsilon_{3},\delta,2\delta,\delta\pm\varepsilon_{1},\pm\delta+\varepsilon_{2},\pm\delta-\varepsilon_{3}\rbrace$ & $\varepsilon_{2}-\varepsilon_{3}=3\alpha_{1}+2\alpha_{2}+2\alpha_{3}$\tabularnewline
\hline 
\end{longtable}
\end{doublespace}

\noindent \textbf{5. $F(4)$} 

The Lie superalgebra $\mathfrak{g}=F(4)$ has the root system $\Phi=\Phi_{\bar{0}}\cup\Phi_{\bar{1}}$
such that $\Phi_{\bar{0}}=\{\pm\delta,\pm\varepsilon_{i}\pm\varepsilon_{j},\pm\varepsilon_{i}:i\neq j,i,j=1,2,3\}$
and $\Phi_{\bar{1}}=\{\frac{1}{2}(\pm\delta\pm\varepsilon_{1}\pm\varepsilon_{2}\pm\varepsilon_{3})\}$,
where $\{\delta,\varepsilon_{1,}\varepsilon_{2},\varepsilon_{3}\}$
is an orthogonal basis of $\mathbb{R}\Phi$ such that $(\delta,\delta)=-6$
and $(\varepsilon_{i},\varepsilon_{j})=2\delta_{ij}$. According to
\cite[Subsection 6.4]{han-exp}, there are six conjugacy classes of
systems of simple roots which we list in the table below.

\begin{doublespace}
\noindent %
\begin{longtable}[c]{|>{\centering}m{4cm}|>{\centering}m{5cm}|>{\centering}m{4cm}|}
\caption{\label{tab:highest root F}The highest root in system of positive
roots for $F(4)$}
\tabularnewline
\endfirsthead
\hline 
System of simple roots $\Pi=\{\alpha_{1},\alpha_{2},\alpha_{3},\alpha_{4}\}$ & Associated system of positive roots $\Phi^{+}$ & The highest root $\widetilde{\alpha}$\tabularnewline
\hline 
\hline 
$\{\frac{1}{2}(\delta-\varepsilon_{1}-\varepsilon_{2}-\varepsilon_{3}),\varepsilon_{3},\varepsilon_{2}-\varepsilon_{3},\varepsilon_{1}-\varepsilon_{2}\}$ & $\{\varepsilon_{1},\varepsilon_{2},\varepsilon_{3},\varepsilon_{1}\pm\varepsilon_{2},\varepsilon_{1}\pm\varepsilon_{3},\varepsilon_{2}\pm\varepsilon_{3},\delta,\frac{1}{2}(\delta\pm\varepsilon_{1}\pm\varepsilon_{2}\pm\varepsilon_{3})\}$ & $\delta=2\alpha_{1}+3\alpha_{2}+2\alpha_{3}+\alpha_{4}$\tabularnewline
\hline 
$\{\frac{1}{2}(-\delta+\varepsilon_{1}+\varepsilon_{2}+\varepsilon_{3}),\frac{1}{2}(\delta-\varepsilon_{1}-\varepsilon_{2}+\varepsilon_{3}),\varepsilon_{2}-\varepsilon_{3},\varepsilon_{1}-\varepsilon_{2}\}$ & $\{\varepsilon_{1},\varepsilon_{2},\varepsilon_{3},\varepsilon_{1}\pm\varepsilon_{2},\varepsilon_{1}\pm\varepsilon_{3},\varepsilon_{2}\pm\varepsilon_{3},\delta,\frac{1}{2}(\pm\delta+\varepsilon_{1}+\varepsilon_{2}+\varepsilon_{3}),\frac{1}{2}(\delta+\varepsilon_{1}\pm\varepsilon_{2}-\varepsilon_{3}),\frac{1}{2}(\delta-\varepsilon_{1}\pm\varepsilon_{2}+\varepsilon_{3}),\frac{1}{2}(\delta-\varepsilon_{1}+\varepsilon_{2}-\varepsilon_{3}),\frac{1}{2}(\delta+\varepsilon_{1}-\varepsilon_{2}+\varepsilon_{3})\}$ & $\frac{1}{2}(\delta+\varepsilon_{1}+\varepsilon_{2}+\varepsilon_{3})=2\alpha_{1}+3\alpha_{2}+2\alpha_{3}+\alpha_{4}$\tabularnewline
\hline 
$\{\varepsilon_{1}-\varepsilon_{2},\frac{1}{2}(\delta-\varepsilon_{1}+\varepsilon_{2}-\varepsilon_{3}),\frac{1}{2}(-\delta+\varepsilon_{1}+\varepsilon_{2}-\varepsilon_{3}),\varepsilon_{3}\}$ & $\{\varepsilon_{1},\varepsilon_{2},\varepsilon_{3},\varepsilon_{1}\pm\varepsilon_{2},\varepsilon_{1}\pm\varepsilon_{3},\varepsilon_{2}\pm\varepsilon_{3},\delta,\frac{1}{2}(\pm\delta+\varepsilon_{1}+\varepsilon_{2}\pm\varepsilon_{3}),\frac{1}{2}(\delta-\varepsilon_{1}+\varepsilon_{2}\pm\varepsilon_{3}),\frac{1}{2}(\delta+\varepsilon_{1}-\varepsilon_{2}\pm\varepsilon_{3})\}$ & $\varepsilon_{1}+\varepsilon_{2}=\alpha_{1}+2\alpha_{2}+2\alpha_{3}+2\alpha_{4}$\tabularnewline
\hline 
$\{\frac{1}{2}(\delta+\varepsilon_{1}-\varepsilon_{2}-\varepsilon_{3}),\frac{1}{2}(\delta-\varepsilon_{1}+\varepsilon_{2}+\varepsilon_{3}),\frac{1}{2}(-\delta+\varepsilon_{1}-\varepsilon_{2}+\varepsilon_{3}),\varepsilon_{2}-\varepsilon_{3}\}$ & $\{\varepsilon_{1},\varepsilon_{2},\varepsilon_{3},\varepsilon_{1}\pm\varepsilon_{2},\varepsilon_{1}\pm\varepsilon_{3},\varepsilon_{2}\pm\varepsilon_{3},\delta,\frac{1}{2}(\pm\delta+\varepsilon_{1}\pm\varepsilon_{2}+\varepsilon_{3}),\frac{1}{2}(\pm\delta+\varepsilon_{1}+\varepsilon_{2}-\varepsilon_{3}),\frac{1}{2}(\delta+\varepsilon_{1}-\varepsilon_{2}-\varepsilon_{3}),\frac{1}{2}(\delta-\varepsilon_{1}+\varepsilon_{2}+\varepsilon_{3})$ & $\varepsilon_{1}+\varepsilon_{2}=\alpha_{1}+2\alpha_{2}+3\alpha_{3}+2\alpha_{4}$\tabularnewline
\hline 
$\{\delta,\frac{1}{2}(-\delta+\varepsilon_{1}-\varepsilon_{2}-\varepsilon_{3}),\varepsilon_{3},\varepsilon_{2}-\varepsilon_{3}\}$ & $\{\varepsilon_{1},\varepsilon_{2},\varepsilon_{3},\varepsilon_{1}\pm\varepsilon_{2},\varepsilon_{1}\pm\varepsilon_{3},\varepsilon_{2}\pm\varepsilon_{3},\delta,\frac{1}{2}(\pm\delta+\varepsilon_{1}\pm\varepsilon_{2}\pm\varepsilon_{3})\}$ & $\varepsilon_{1}+\varepsilon_{2}=\alpha_{1}+2\alpha_{2}+3\alpha_{3}+2\alpha_{4}$\tabularnewline
\hline 
$\{\delta,\frac{1}{2}(-\delta-\varepsilon_{1}+\varepsilon_{2}+\varepsilon_{3}),\varepsilon_{1}-\varepsilon_{2},\varepsilon_{2}-\varepsilon_{3}\}$ & $\{\varepsilon_{1},\varepsilon_{2},\varepsilon_{3},\varepsilon_{1}\pm\varepsilon_{2},\varepsilon_{1}\pm\varepsilon_{3},\varepsilon_{2}\pm\varepsilon_{3},\delta,\frac{1}{2}(\pm\delta+\varepsilon_{1}+\varepsilon_{2}\pm\varepsilon_{3}),\frac{1}{2}(\pm\delta-\varepsilon_{1}+\varepsilon_{2}+\varepsilon_{3}),\frac{1}{2}(\pm\delta+\varepsilon_{1}-\varepsilon_{2}+\varepsilon_{3})\}$  & $\varepsilon_{1}+\varepsilon_{2}=2\alpha_{1}+4\alpha_{2}+3\alpha_{3}+2\alpha_{4}$\tabularnewline
\hline 
\end{longtable}
\end{doublespace}

From above we know that the highest root depends on the choice of
simple roots. In this paper we define a good prime $p$ as below.
\begin{defn}
\label{def:good prime}A prime $p$ is said to be \textit{good} for
$\mathfrak{g}$ if it is greater than all coefficients $a_{i}$ of
$\alpha_{i}$ in the highest root of every system of positive roots
$\Phi^{+}$, and to be \textit{bad} for $\mathfrak{g}$ otherwise. 

Thus if $\mathfrak{g}$ is of type $A(m,n)$, then all primes are
good; if $\mathfrak{g}$ is of type $B(m,n)$, $C(n)$, $D(m,n)$
or $D(2,1;\alpha)$, then $2$ is bad; if $\mathfrak{g}=G(3)$ or
$F(4)$, then $2$ and $3$ are bad.
\end{defn}

\subsection{Cocharacters associated to nilpotent elements}

Let $G_{\bar{0}}$ be a reductive algebraic group over $\mathbb{K}$
defined as in Table \ref{tab:algebraic-groups} so that the adjoint
action of $G_{\bar{0}}$ on $\mathfrak{g}$. Recall that any homomorphism
of algebraic groups $\tau:\mathbb{K}^{\times}\rightarrow G_{\bar{0}}$
defines a grading $\mathfrak{g}=\bigoplus_{j\in\mathbb{Z}}\mathfrak{g}(j;\tau)$
such that 
\begin{equation}
\mathfrak{g}(j;\tau)=\{x\in\mathfrak{g}:\mathrm{Ad}(\tau(t))(x)=t^{j}x\text{ for all }t\in\mathbb{K}^{\times}\}.\label{eq:cocharacter}
\end{equation}
 A \textit{parabolic subgroup} $P$ of $G_{\bar{0}}$ is a subgroup
containing a maximal connected solvable algebraic subgroup of $G_{\bar{0}}$.
By \cite[Theorem 30.2]{Humphreys2012}, any parabolic subgroup $P$
of $G_{\bar{0}}$ has a Levi decomposition $P=LV$ where $L$ is a
\textit{Levi factor}. Throughout this paper we call a Levi factor
of a parabolic subgroup of $G_{\bar{0}}$ a \textit{Levi subgroup}
\textit{of} $G_{\bar{0}}$. Let $L$ be a Levi subgroup of $G_{\bar{0}}$.
A nilpotent element $e\in\mathfrak{g}_{\bar{0}}$ is called \textit{distinguished}
in $\mathrm{Lie}(L)$ if each torus contained in $L^{e}$ is also
contained in the centre $Z(L)$ of $L$. According to \cite[Section 2]{Lawther2008},
every nilpotent element in $\mathfrak{g}_{\bar{0}}$ is distinguished
in the Lie algebra of some Levi subgroup of $G_{\bar{0}}$. 
\begin{defn}
\cite[Definition 5.3]{Janzten} \label{def:cocharacter}Let $e\in\mathfrak{g}_{\bar{0}}$
be nilpotent. A cocharacter $\tau:\mathbb{K}^{\times}\rightarrow G_{\bar{0}}$
is called \textit{associated to }$e$ if $e\in\mathfrak{g}(2;\tau)$
and there exists a Levi subgroup $L$ of $G_{\bar{0}}$ such that
$e$ is distinguished in $L$ and $\mathrm{im}(\tau)\subseteq[L,L]$. 
\end{defn}

By \cite[Lemma 5.3]{Janzten}, cocharacters associated to $e$ do
exist if $\mathrm{char}(\mathbb{K})$ is good for $G_{\bar{0}}$,
and two cocharacters associated to $e$ are conjugate under the identity
component $(G_{\bar{0}}^{e})^{\circ}$ of $G_{\bar{0}}^{e}$. 
\begin{thm}
\label{thm:g^e=00003Ddimg(0)+dimg(1)}\cite[Theorem 3.1]{WANGZHAO}
Let $\mathfrak{g}=\mathfrak{g}_{\bar{0}}\oplus\mathfrak{g}_{\bar{1}}$
be one of the basic classical Lie superalgebras in Table \ref{tab:Basic-classical}
with $p>3$ if $\mathfrak{g}=D(2,1;\alpha)$ and $p>15$ if $\mathfrak{g}=G(3)$
or $F(4)$. Let $e\in\mathfrak{g}_{\bar{0}}$ be nilpotent. Then the
cocharacter $\tau$ associated to $e$ defines a $\mathbb{Z}$-grading
$\mathfrak{g}=\bigoplus_{j\in\mathbb{Z}}\mathfrak{g}(j;\tau)$ such
that 
\begin{equation}
e\in\mathfrak{g}(2;\tau)\text{ and }\mathfrak{g}^{e}=\bigoplus_{j\in\mathbb{Z}}\mathfrak{g}^{e}(j;\tau)\text{ where }\mathfrak{g}^{e}(j;\tau)=\mathfrak{g}^{e}\cap\mathfrak{g}(j;\tau);\label{eq:thm 2-1}
\end{equation}
\begin{equation}
\mathfrak{g}^{e}(j;\tau)=0\text{ for all }j<0\text{ and }\dim\mathfrak{g}^{e}=\dim\mathfrak{g}(0;\tau)+\dim\mathfrak{g}(1;\tau).\label{eq:thm 2-2}
\end{equation}
\end{thm}

Given a nilpotent element $e\in\mathfrak{g}_{\bar{0}}$, we fix the
associated cocharacter $\tau$ as follows. Let $\mathfrak{h}$ be
a Cartan subalgebra of $\mathfrak{g}$. Recall that the non-degenerate
supersymmetric invariant even bilinear form $(\cdotp,\cdotp)$ on
$\mathfrak{g}$ restricts to a non-degenerate symmetric bilinear form
on $\mathfrak{h}$. This allows us to identify $\mathfrak{h}\cong\mathfrak{h}^{*}$
and obtain a symmetric bilinear form $(\cdotp,\cdotp)$ on $\mathfrak{h}^{*}$.
For $\alpha,\beta\in\mathfrak{h}^{*}$, we define 
\begin{equation}
\langle\beta,\alpha\rangle=\begin{cases}
\frac{2(\alpha,\beta)}{(\alpha,\alpha)} & \text{if }(\alpha,\alpha)\neq0;\\
-1 & \text{if }(\alpha,\alpha)=0\text{ and }(\alpha,\beta)\neq0;\\
0 & \text{if }(\alpha,\alpha)=(\alpha,\beta)=0
\end{cases}\label{eq:<beta,alpha>}
\end{equation}
 as in \cite[Section 2.3]{Frappat1989}. Fix a maximal torus $T$
of $G_{\bar{0}}$. Then we consider a cocharacter $\tau:\mathbb{K}^{\times}\rightarrow G_{\bar{0}}$
associated to $e$ in a similar way to \cite[Section 6]{Lawther2008}
such that $\tau\left(t\right)=\prod_{i=1}^{s}h_{\alpha_{i}}(t^{c_{i}})$
where $c_{i}$ is defined in Section \ref{sec:Introduction}. Note
that $h_{\alpha_{i}}(t)\in T$ is determined by $h_{\alpha_{i}}(t)\cdotp e_{\alpha}=t^{\langle\alpha,\alpha_{i}\rangle}e_{\alpha}$
for all $t\in\mathbb{K}^{\times}$. 
\begin{lem}
\label{lem:g(j)}Let $\alpha\in\Phi$. If $e_{\alpha}\in\mathfrak{g}_{\mathbb{Z}}(j;\mathrm{ad}h)$,
then $e_{\alpha}\in\mathfrak{g}(j;\tau)$.
\end{lem}

\begin{proof}
Suppose $e_{\alpha}\in\mathfrak{g}_{\mathbb{Z}}(j;\mathrm{ad}h)$,
we know that 
\[
[h,e_{\alpha}]=[\sum_{i=1}^{s}c_{i}h_{i},e_{\alpha}]=\sum_{i=1}^{s}c_{i}[h_{i},e_{\alpha}]=\sum_{i=1}^{s}c_{i}\langle\alpha,\alpha_{i}\rangle e_{\alpha}.
\]
 Hence, we have that $\sum_{i=1}^{s}c_{i}\langle\alpha,\alpha_{i}\rangle=j$.
We also have that 
\[
\tau\left(t\right)\cdotp e_{\alpha}=\prod_{i=1}^{s}h_{\alpha_{i}}(t^{c_{i}})\cdotp e_{\alpha}=\prod_{i=1}^{s}t^{\langle\alpha,\alpha_{i}\rangle}e_{\alpha}=t^{\sum_{i=1}^{s}\langle\alpha,\alpha_{i}\rangle}e_{\alpha}\in\mathfrak{g}(j;\tau).
\]
 Therefore, we obtain that $e_{\alpha}\in\mathfrak{g}(j;\tau)$, as
required.
\end{proof}
From now on let us write $\mathfrak{g}(j)$ for $\mathfrak{g}(j;\tau)$
and $\mathfrak{g}^{e}(j)$ for $\mathfrak{g}^{e}(j;\tau)$. We also
denote by $\mathfrak{g}_{\bar{i}}^{e}(j)=\mathfrak{g}_{\bar{i}}\cap\mathfrak{g}^{e}(j)$.

\section{Centralizer and Centre of centralizer of nilpotent elements in Lie
superalgebras of type $A(m,n)$\label{sec:A(m,n)}}

\subsection{Centralizer of nilpotent even elements in Lie superalgebras $A(m,n)$}

In this subsection, we recall the centralizer of a nilpotent even
element $e$ in a Lie superalgebra of type $A(m,n)$ which is analogous
to the Lie algebra case in \cite[Section 3]{Janzten}. For $\mathfrak{gl}(m|n)$
this was done in \cite[Subsection 3.2]{Hoyt2012} for a field of zero
characteristic and the construction is identical for a field of prime
characteristic by \cite[Subsection 3.2]{WANGZHAO}.

\begin{singlespace}
Let $V=V_{\bar{0}}\oplus V_{\bar{1}}$ be a finite-dimensional vector
superspace over $\mathbb{K}$ such that $\dim V_{\bar{0}}=m$ and
$\dim V_{\bar{1}}=n$. Let $\bar{\mathfrak{g}}=\bar{\mathfrak{g}}_{\bar{0}}\oplus\bar{\mathfrak{g}}_{\bar{1}}=\mathfrak{gl}(V)=\mathfrak{gl}(m|n)$
and $e\in\bar{\mathfrak{g}}_{\bar{0}}$ be nilpotent. Note that the
nilpotent orbits in $\bar{\mathfrak{g}}_{\bar{0}}$ are parametrized
by the partitions of $(m|n)$. Let $\lambda$ be a partition of $(m|n)$
such that 
\begin{equation}
\lambda=(p_{1},\dots,p_{r}\mid q_{1},\dots,q_{s})\label{eq:gl(m,n)-partition 0}
\end{equation}
where $p_{1}\geq\dots\geq p_{r},q_{1}\geq\dots\geq q_{s}$, $\sum_{i=1}^{r}p_{i}=m$
and $\sum_{i=1}^{s}q_{i}=n$. By rearranging the order of numbers
in (\ref{eq:gl(m,n)-partition 0}), we write 
\begin{equation}
\lambda=(\lambda_{1},\dots,\lambda_{r+s})\label{eq:gl(m,n)-partition 1}
\end{equation}
 such that $\lambda_{1}\geq\dots\geq\lambda_{r+s}$. For $i=1,\dots,r+s$,
we define $\bar{i}\in\{\bar{0},\bar{1}\}$ such that for $c\in\mathbb{Z}$,
we have that $\ensuremath{\left|\{i:\lambda_{i}=c,\bar{i}=\bar{0}\}\right|}=\ensuremath{\left|\{j:p_{j}=c\}\right|}$,
$\left|\{i:\lambda_{i}=c,\bar{i}=\bar{1}\}\right|=\ensuremath{\left|\{j:q_{j}=c\}\right|}$
and if $\lambda_{i}=\lambda_{j}$, $\bar{i}=\bar{0},\bar{j}=\bar{1}$,
then $i<j$. i.e. $\sum_{\ensuremath{\bar{i}}=\bar{0}}\lambda_{i}=m$
and $\sum_{\bar{i}=\bar{1}}\lambda_{i}=n$.
\end{singlespace}

Recall that $\bar{\mathfrak{g}}=\mathrm{End}(V_{\bar{0}}\oplus V_{\bar{1}})$
where $\bar{\mathfrak{g}}_{\bar{0}}=\mathrm{End}(V_{\bar{0}})\oplus\mathrm{End}(V_{\bar{1}})$
and $\bar{\mathfrak{g}}_{\bar{1}}=\mathrm{Hom}(V_{\bar{0}},V_{\bar{1}})\oplus\mathrm{Hom}(V_{\bar{1}},V_{\bar{0}})$.
Then there exist $v_{1},\dots,v_{r+s}\in V$ such that $V_{\bar{0}}=\langle e^{k}v_{i}:0\leq k\leq\lambda_{i}-1,\bar{i}=\bar{0}\rangle$,
$V_{\bar{1}}=\langle e^{k}v_{i}:0\leq k\leq\lambda_{i}-1,\bar{i}=\bar{1}\rangle$
and $e^{\lambda_{i}}v_{i}=0$ for $1\leq i\leq r+s$ by \cite[Section 3]{Janzten}.
According to \cite[Subsection 3.2]{WANGZHAO}, the following elements
\begin{equation}
\left\{ \xi_{i}^{j,k}:1\leq i,j\leq r+s\text{ and }\max\{\lambda_{j}-\lambda_{i},0\}\leq k\leq\lambda_{j}-1\right\} \label{eq:basis gl}
\end{equation}
 form a basis for $\bar{\mathfrak{g}}^{e}$ where $\xi_{i}^{j,k}(v_{t})=\delta_{it}e^{k}v_{j}$.
\begin{example}
Suppose $\bar{\mathfrak{g}}=\bar{\mathfrak{g}}_{\bar{0}}\oplus\bar{\mathfrak{g}}_{\bar{1}}=\mathfrak{gl}(7|3)=\mathfrak{gl}(V)$.
Let $e\in\bar{\mathfrak{g}}_{\bar{0}}$ be a nilpotent element that
is parametrized by the partition $\lambda=(5,2|3)$. Rewrite $\lambda$
using the notation defined in (\ref{eq:gl(m,n)-partition 1}), we
have that $\lambda_{1}=5,\lambda_{2}=3$ and $\lambda_{3}=2$. Denote
by $J_{\lambda_{i}}$ the $(\lambda_{i}\times\lambda_{i})$-matrix
where the $(j,j+1)$ entries with $1\leq j<\lambda_{i}$ are equal
to $1$ and all remaining entries are equal to zero. Then the order
of basis for $V$ is rearranged to be $(e^{4}v_{1},e^{3}v_{1},e^{2}v_{1},ev_{1},v_{1},e^{2}v_{2},ev_{2},v_{2},ev_{3},v_{3})$
such that $v_{1},v_{3}\in V_{\bar{0}}$, $v_{2}\in V_{\bar{1}}$ and
with respect to this basis $e$ has block form 
\[
\begin{pmatrix}J_{\lambda_{1}} & 0 & 0\\
0 & J_{\lambda_{2}} & 0\\
0 & 0 & J_{\lambda_{3}}
\end{pmatrix}.
\]
 In Figure \ref{fig:Example}, we write elements $\xi_{i}^{j,k}\in\bar{\mathfrak{g}}^{e}$
in block matrices with respect to the basis for $V$. 

\begin{figure}[H]
\includegraphics[scale=0.35]{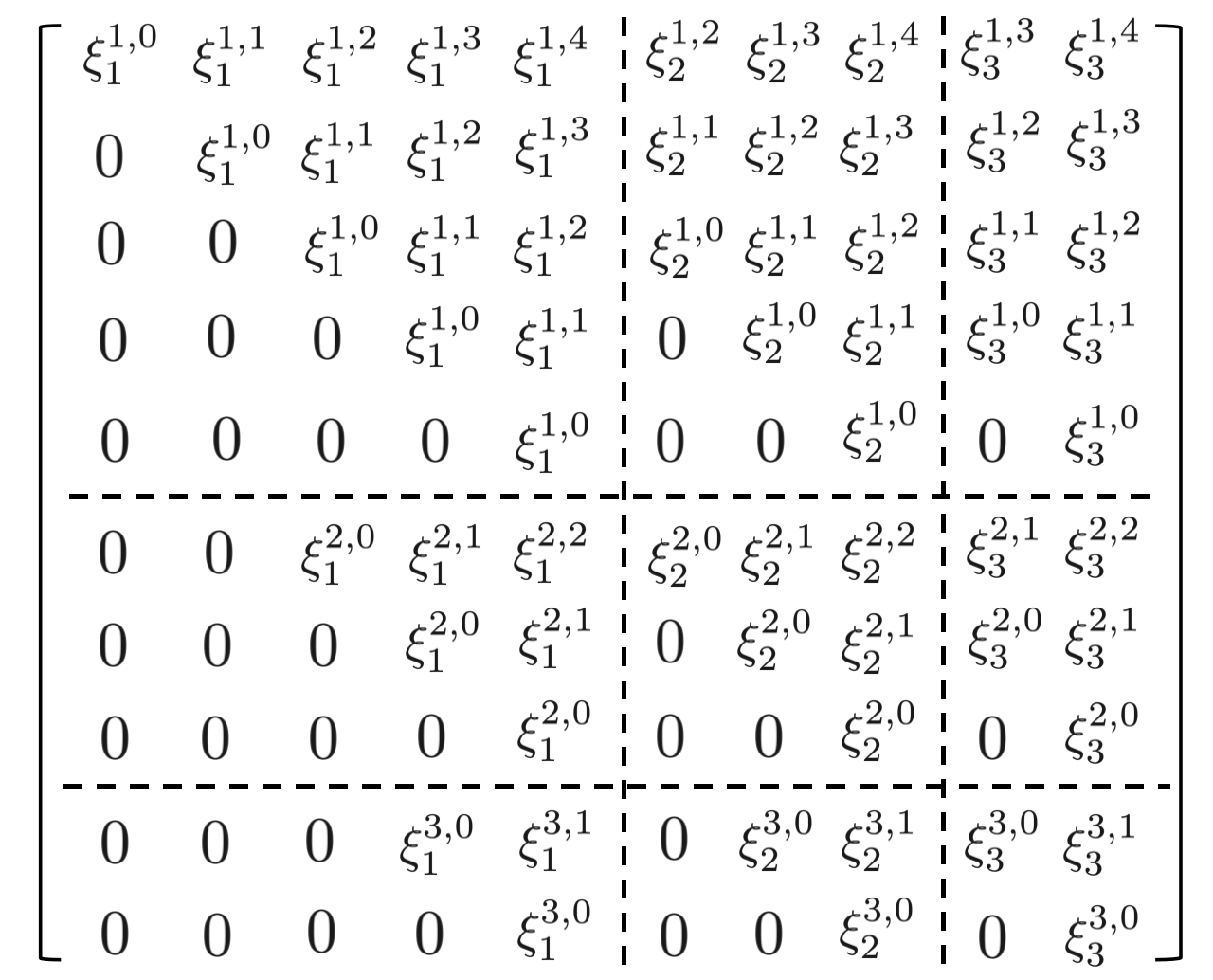}

\caption{\label{fig:Example}$\xi_{i}^{j,k}$}
\end{figure}
\end{example}

Now let $\mathfrak{g}=\mathfrak{g}_{\bar{0}}\oplus\mathfrak{g}_{\bar{1}}=\mathfrak{sl}(V)=\mathfrak{sl}(m|n)$.
Based on the above example, in order to determine $\mathfrak{g}^{e}$,
we need to let $\sum_{i=1}^{r+s}(-1)^{\bar{i}}\lambda_{i}\xi_{i}^{i,0}=0$.
Note that if $\mathrm{char}(\mathbb{K})=p$ divides $\lambda_{i}$
for all $1\leq i\leq r+s$, we have that $\sum_{i=1}^{r+s}(-1)^{\bar{i}}\lambda_{i}\xi_{i}^{i,0}=0$.
Hence, we deduce that $\dim\mathfrak{g}^{e}=\dim\bar{\mathfrak{g}}^{e}$
and the construction of $\bar{\mathfrak{g}}^{e}$ is identical to
$\mathfrak{g}^{e}$ for this case. However, if there exist some $\lambda_{i}$
which is not divisible by $p$, we have that $\dim\mathfrak{g}^{e}=\dim\bar{\mathfrak{g}}^{e}-1$
as $\sum_{i=1}^{r+s}(-1)^{\bar{i}}\lambda_{i}\xi_{i}^{i,0}$ has to
be zero. Thus the diagonal matrices $\xi_{i}^{i,0}$ are not in $\mathfrak{g}^{e}$
and a basis for $\mathfrak{g}^{e}$ is of the form 
\[
\left\{ \xi_{i}^{j,k}:1\leq i\neq j\leq r+s\text{ and }\max\{\lambda_{j}-\lambda_{i},0\}\leq k\leq\lambda_{j}-1\right\} \cup\left\{ \xi_{i}^{i,k}:0<k\leq\lambda_{i}-1\right\} 
\]
\[
\cup\left\{ \frac{\lambda_{i+1}}{p^{a}}\xi_{i}^{i,0}-(-1)^{\bar{i}}\frac{\lambda_{i}}{p^{a}}\xi_{i+1}^{i+1,0}:1\leq i<r+s,a\text{ is maximal such that }p^{a}\text{ divides both }\lambda_{i}\text{ and }\lambda_{i+1}\right\} .
\]

Note that when $m=n$, we consider $\mathfrak{g}=\mathfrak{g}_{\bar{0}}\oplus\mathfrak{g}_{\bar{1}}=\mathfrak{psl}(n|n)$
instead of $\mathfrak{sl}(n|n)$. Let $e\in\mathfrak{g}_{\bar{0}}$
be nilpotent. By abuse of notation, we use the same letter $e$ for
its lift in $\mathfrak{sl}(n|n)$.
\begin{lem}
When $p$ does not divide $n$, we have that $\dim\mathfrak{g}^{e}=\dim\mathfrak{sl}(n|n)^{e}-1$.
\end{lem}

\begin{proof}
Let $\phi:\mathfrak{sl}(n|n)\rightarrow\mathfrak{g}$ be the quotient
map. For $x=x_{\bar{0}}+x_{\bar{1}}$ such that $x+\mathfrak{z}(\mathfrak{sl}(n|n))\in\mathfrak{g}^{e}$,
we have that $[e,x]\in\mathfrak{z}(\mathfrak{sl}(n|n))$. Since $[e,x]=[e,x_{\bar{0}}]+[e,x_{\bar{1}}]$,
then $[e,x]\in\mathfrak{z}(\mathfrak{sl}(n|n))$ implies that $[e,x_{\bar{1}}]=0$
and $[e,x_{\bar{0}}]\in\mathfrak{z}(\mathfrak{sl}(n|n))$. When $p$
does not divide $n$, this is impossible since $[e,x_{\bar{0}}]\in[\mathfrak{g}_{\bar{0}},\mathfrak{g}_{\bar{0}}]=\mathfrak{sl}_{n}(\mathbb{K})\oplus\mathfrak{sl}_{n}(\mathbb{K})$
and $\mathfrak{z}(\mathfrak{sl}(n|n))\notin\mathfrak{sl}_{n}(\mathbb{K})\oplus\mathfrak{sl}_{n}(\mathbb{K})$.
Hence, we have that $\phi^{e}:\mathfrak{sl}(n|n)^{e}\rightarrow\mathfrak{g}^{e}$
is surjective and $\mathrm{ker}(\phi^{e})=\mathfrak{z}(\mathfrak{sl}(n|n))$.
Therefore, we deduce that $\dim\mathfrak{g}^{e}=\dim\mathfrak{sl}(n|n)^{e}-1$. 
\end{proof}

Since $\mathfrak{z}(\mathfrak{sl}(n|n))\subseteq\mathfrak{sl}(n|n)^{e}$
and $\phi$ sends $\mathfrak{z}(\mathfrak{sl}(n|n))$ to zero, we
can write a basis for $\mathfrak{g}^{e}$ as 
\begin{equation}
\left\{ \xi_{i}^{j,k}:1\leq i\neq j\leq r+s\text{ and }\max\{\lambda_{j}-\lambda_{i},0\}\leq k\leq\lambda_{j}-1\right\} \label{eq:psl(n|n)^e}
\end{equation}
\[
\cup\left\{ \xi_{i}^{i,k}:1\leq i\leq r+s,0<k\leq\lambda_{i}-1\right\} 
\]
\[
\cup\left\{ \frac{\lambda_{i+1}}{p^{a}}\xi_{i}^{i,0}-(-1)^{\bar{i}}\frac{\lambda_{i}}{p^{a}}\xi_{i+1}^{i+1,0}:1<i<r+s,a\text{ is maximal such that }p^{a}\text{ divides both }\lambda_{i}\text{ and }\lambda_{i+1}\right\} .
\]

\begin{rem}
Note that when $p$ divides $n$, then $\mathfrak{z}(\mathfrak{sl}(n|n))\subseteq\mathfrak{sl}_{n}(\mathbb{K})\oplus\mathfrak{sl}_{n}(\mathbb{K})$
and $\phi^{e}:\mathfrak{sl}(n|n)^{e}\rightarrow\mathfrak{g}^{e}$
is not surjective. In this paper, we rule out this case when considering
$\mathfrak{psl}(n|n)^{e}$ and $\mathfrak{z}(\mathfrak{psl}(n|n)^{e})$.
\end{rem}

\subsection{Centre of centralizer of nilpotent even elements in Lie superalgebras
$A(m,n)$}

The centre of centralizer of a nilpotent even element $e$ in $\mathfrak{gl}(m|n)$
was done for a field of zero characteristic in \cite[Subsection 3.3]{han-classical},
we prove that the structure is identical over our field $\mathbb{K}$
in the following theorem.
\begin{thm}
\label{thm:centre gl(m|n)}Let $\bar{\mathfrak{g}}=\bar{\mathfrak{g}}_{\bar{0}}\oplus\bar{\mathfrak{g}}_{\bar{1}}=\mathfrak{gl}(m|n)$.
Let $e\in\bar{\mathfrak{g}}_{\bar{0}}$ be nilpotent and $\lambda=(\lambda_{1},\dots,\lambda_{r+s})$
be the partition with respect to $e$ defined as in \ref{eq:gl(m,n)-partition 1}.
Then $\mathfrak{z}(\bar{\mathfrak{g}}^{e})=\langle I,e,\dots,e^{\lambda_{1}-1}\rangle$.
\end{thm}

\begin{proof}
Let us denote by $I_{m|0}\in\bar{\mathfrak{g}}^{e}$ be the $(m+n)\times(m+n)$
matrix such that 
\[
I_{m|0}=\begin{pmatrix}I_{m} & 0_{m\times n}\\
0_{n\times m} & 0_{n\times n}
\end{pmatrix}
\]
 where $I_{m}$ is the $m\times m$ identity matrix and $0_{m\times n},0_{n\times m},0_{n\times n}$
are matrices with zero entries. For any $x\in\mathfrak{z}(\bar{\mathfrak{g}}^{e})$,
we have that $[I_{m|0},x]=0$. Hence, we obtain that $\mathfrak{z}(\bar{\mathfrak{g}}^{e})\subseteq\bar{\mathfrak{g}}^{I_{m|0}}$.
Let $y=\begin{pmatrix}A & B\\
C & D
\end{pmatrix}\in\bar{\mathfrak{g}}$, we compute that $[I_{m|0},y]=0$ if and only if $B=C=0$. Thus we
deduce that $y\in\bar{\mathfrak{g}}_{\bar{0}}$ and $\bar{\mathfrak{g}}^{I_{m|0}}=\bar{\mathfrak{g}}_{\bar{0}}$.
Therefore, we have that $\mathfrak{z}(\bar{\mathfrak{g}}^{e})\subseteq\bar{\mathfrak{g}}_{\bar{0}}$.

Observe that $\bar{\mathfrak{g}}$ can be viewed as a vector space
that is isomorphic to $\mathfrak{gl}_{m+n}(\mathbb{K})$. We also
have that $\bar{\mathfrak{g}}^{e}$ is isomorphic to $\mathfrak{gl}_{m+n}(\mathbb{K})^{e}$
as a vector space because $[e,x]=[e,x']$ for $x\in\bar{\mathfrak{g}}$
and $x'\in\mathfrak{gl}_{m+n}(\mathbb{K})$. Let us denote 
\[
\mathfrak{z}(\mathfrak{gl}_{m+n}(\mathbb{K})^{e})_{0}=\mathfrak{z}(\mathfrak{gl}_{m+n}(\mathbb{K})^{e})\cap\left\{ \begin{pmatrix}A' & 0\\
0 & D'
\end{pmatrix}\in\mathfrak{gl}_{m+n}(\mathbb{K}):A'\in\mathfrak{gl}_{m}(\mathbb{K})\text{ and }D'\in\mathfrak{gl}_{n}(\mathbb{K})\right\} ,
\]
\[
\text{and }\mathfrak{z}(\mathfrak{gl}_{m+n}(\mathbb{K})^{e})_{1}=\mathfrak{z}(\mathfrak{gl}_{m+n}(\mathbb{K})^{e})\cap\left\{ \begin{pmatrix}0 & B'\\
C' & 0
\end{pmatrix}:B'\text{ is a }m\times n\text{ matrix, }C'\text{ is a }n\times m\text{ matrix}\right\} .
\]
 Let $\mathfrak{z}(\bar{\mathfrak{g}}^{e})=\mathfrak{z}(\bar{\mathfrak{g}}^{e})_{\bar{0}}\oplus\mathfrak{z}(\bar{\mathfrak{g}}^{e})_{\bar{1}}.$
Note that $\mathfrak{z}(\bar{\mathfrak{g}}^{e})_{\bar{0}}=\mathfrak{z}(\mathfrak{gl}_{m+n}(\mathbb{K})^{e})_{0}$
by definition. Since $\mathfrak{z}(\bar{\mathfrak{g}}^{e})\subseteq\bar{\mathfrak{g}}_{\bar{0}}$,
we have that $\mathfrak{z}(\bar{\mathfrak{g}}^{e})_{\bar{1}}=0$.
Using a similar argument as in the first part of this proof, we also
obtain that $\mathfrak{z}(\mathfrak{gl}_{m+n}(\mathbb{K})^{e})_{1}=0$.
By \cite[Section 4]{Lawther2008} we know that $\mathfrak{z}(\mathfrak{gl}_{m+n}(\mathbb{K})^{e})=\langle I,e,\dots,e^{\lambda_{1}-1}\rangle$.
Therefore, we obtain that $\mathfrak{z}(\bar{\mathfrak{g}}^{e})=\mathfrak{z}(\bar{\mathfrak{g}}^{e})_{\bar{0}}=\mathfrak{z}(\mathfrak{gl}_{m+n}(\mathbb{K})^{e})=\mathfrak{z}(\mathfrak{gl}_{m+n}(\mathbb{K})^{e})_{0}=\langle I,e,\dots,e^{\lambda_{1}-1}\rangle$.
\end{proof}
Next we look at the centre of centralizer of a nilpotent element $e\in\mathfrak{sl}(m|n)_{\bar{0}}$
when $m\neq n$ and if $m=n$, we consider $\mathfrak{psl}(n|n)$
instead.
\begin{thm}
Let $\mathfrak{g}=\mathfrak{g}_{\bar{0}}\oplus\mathfrak{g}_{\bar{1}}=\mathfrak{sl}(m|n)$
for $m\neq n$. Let $e\in\mathfrak{g}_{\bar{0}}$ be nilpotent and
$\lambda=(\lambda_{1},\dots,\lambda_{r+s})$ be the partition with
respect to $e$ defined as in \ref{eq:gl(m,n)-partition 1}. If $\mathrm{char}(\mathbb{K})=p$
divides $\lambda_{i}$ for all $1\leq i\leq r+s$, then $\mathfrak{z}(\mathfrak{g}^{e})=\langle I,e,e^{2},\dots,e^{\lambda_{1}-1}\rangle$;
If there exists some $\lambda_{i}$ which is not divisible by $p$,
then $\mathfrak{z}(\mathfrak{g}^{e})=\langle e,e^{2},\dots,e^{\lambda_{1}-1}\rangle$.
\end{thm}

\begin{proof}
When $m\neq n$. If $\mathrm{char}(\mathbb{K})=p$ divides $\lambda_{i}$
for all $1\leq i\leq r+s$, then $\mathfrak{g}^{e}=\bar{\mathfrak{g}}^{e}$
and thus by Theorem \ref{thm:centre gl(m|n)} we have that $\mathfrak{z}(\mathfrak{g}^{e})=\langle I,e,e^{2},\dots,e^{\lambda_{1}-1}\rangle$.
If there exists some $\lambda_{i}$ which is not divisible by $p$,
then $\bar{\mathfrak{g}}^{e}\subseteq\mathfrak{g}^{e}\oplus\mathbb{K}I$
and thus $\mathfrak{z}(\bar{\mathfrak{g}}^{e})\subseteq\mathfrak{z}(\mathfrak{g}^{e})\oplus\mathbb{K}I$.
For an element $x\in\mathfrak{z}(\mathfrak{g}^{e})$, observe that
$[x,I]=0$ and $[x,y]=0$ for all $y\in\mathfrak{g}^{e}$. Hence,
we deduce that $[x,y]=0$ for all $y\in\bar{\mathfrak{g}}^{e}$ and
thus $x\in\mathfrak{z}(\bar{\mathfrak{g}}^{e})$. This implies that
$\mathfrak{z}(\mathfrak{g}^{e})\subseteq\mathfrak{z}(\bar{\mathfrak{g}}^{e})$.
Therefore, we obtain that a basis for $\mathfrak{z}(\mathfrak{g}^{e})$
contains all basis vectors of $\mathfrak{z}(\bar{\mathfrak{g}}^{e})$
except the identity matrix $I$.
\end{proof}
In the remaining part of this subsection, let $\mathfrak{g}=\mathfrak{g}_{\bar{0}}\oplus\mathfrak{g}_{\bar{1}}=\mathfrak{psl}(n|n)$
and let $e\in\mathfrak{g}_{\bar{0}}$ be nilpotent. Below we give
a basis of $\mathfrak{z}(\mathfrak{g}^{e})$. Note that a basis of
$\mathfrak{z}(\mathfrak{g}_{\mathbb{C}}^{e})$ is given in \cite[Subsection 4.2]{han-reachable}
and we apply a similar argument in the following theorem. By abuse
of notation, we use the same letter $e$ for its lift in $\mathfrak{sl}(n|n)$.
\begin{thm}
\label{thm:z(g^e)}Let $\mathfrak{g}=\mathfrak{g}_{\bar{0}}\oplus\mathfrak{g}_{\bar{1}}$
be a Lie superalgebra of type $\mathfrak{psl}(n|n)$ over $\mathbb{K}$
such that $\mathrm{char}(\mathbb{K})=p$ does not divide $n$. Let
$\lambda=(\lambda_{1},\dots,\lambda_{r+s})$ be a partition of $(n|n)$
defined as in (\ref{eq:gl(m,n)-partition 1}) and suppose $e\in\mathfrak{g}_{\bar{0}}$
is a nilpotent element with Jordan type $\lambda$. Then $\mathfrak{z}(\mathfrak{g}^{e})=\langle e,\dots,e^{\lambda_{1}-1}\rangle$.
\end{thm}

\begin{proof}
Let us denote by $\mathfrak{h}$ the Cartan subalgebra of $\mathfrak{g}$.
Then elements $\frac{\lambda_{j^{'}}}{p^{a}}\xi_{i^{'}}^{i^{'},0}-(-1)^{\bar{i^{'}}}\frac{\lambda_{i^{'}}}{p^{a}}\xi_{j^{'}}^{j^{'},0}\in\mathfrak{h}$
for $i^{'}\neq j^{'}$ where $a$ is maximal such that $p^{a}$ divides
both $\lambda_{i^{'}}$ and $\lambda_{j^{'}}$. From (\ref{eq:psl(n|n)^e}),
we know that an element $x\in\mathfrak{g}^{e}$ is of the form 
\[
x=\sum_{\substack{1\leq i\neq j\leq r+s\\
\max\{\lambda_{j}-\lambda_{i},0\}\leq k\leq\lambda_{j}-1
}
}a_{i,j,k}\xi_{i}^{j,k}+\sum_{\substack{0<k\leq\lambda_{i}-1\\
1\leq i\leq r+s
}
}b_{i,k}\xi_{i}^{i,k}+\sum_{1<i<r+s}c_{i,0}\left(\frac{\lambda_{i+1}}{p^{a}}\xi_{i}^{i,0}-(-1)^{\bar{i}}\frac{\lambda_{i}}{p^{a}}\xi_{i+1}^{i+1,0}\right)
\]

where $a_{i,j,k},b_{i,k},c_{i,0}\in\mathbb{K}$. 

Now assume $x\in\mathfrak{z}(\mathfrak{g}^{e})$. We divide our analysis
into two cases.

1. When $r+s\geq3$, we know that $[x,\frac{\lambda_{j^{'}}}{p^{a}}\xi_{i^{'}}^{i^{'},0}-(-1)^{\bar{i^{'}}}\frac{\lambda_{i^{'}}}{p^{a}}\xi_{j^{'}}^{j^{'},0}]\in\left\langle \sum_{1\leq i\leq r+s}\xi_{i}^{i,0}\right\rangle $.
We have that 
\begin{align}
[x,\frac{\lambda_{j^{'}}}{p^{a}}\xi_{i^{'}}^{i^{'},0}-(-1)^{\bar{i^{'}}}\frac{\lambda_{i^{'}}}{p^{a}}\xi_{j^{'}}^{j^{'},0}] & =\frac{\lambda_{j^{'}}}{p^{a}}\left(\sum_{j,k}a_{i^{'},j,k}\xi_{i^{'}}^{j,k}-\sum_{i,k}a_{i,i^{'},k}\xi_{i}^{i^{'},k}\right)\nonumber \\
 & -(-1)^{\bar{i^{'}}}\frac{\lambda_{i^{'}}}{p^{a}}\left(\sum_{j,k}a_{j^{'},j,k}\xi_{j^{'}}^{j,k}-\sum_{i,k}a_{i,j^{'},k}\xi_{i}^{j^{'},k}\right).\label{eq:centre psl-1}
\end{align}
 In order to obtain a summand of $\xi_{i}^{i,0}$ with nonzero coefficients
for any $i$ from (\ref{eq:centre psl-1}), we need $k=0$, $j=i^{'}$
(resp. $k=0$, $i=i^{'}$) in the first (resp. second) summand and
we have $a_{i^{'},i^{'},0}\xi_{i^{'}}^{i^{'},0}-a_{i^{'},i^{'},0}\xi_{i^{'}}^{i^{'},0}=0$.
Similarly we have $a_{j^{'},j^{'},0}\xi_{j^{'}}^{j^{'},0}-a_{j^{'},j^{'},0}\xi_{j^{'}}^{j^{'},0}=0$
in the third and fourth summands. Hence, we deduce that $[x,\frac{\lambda_{j^{'}}}{p^{a}}\xi_{i^{'}}^{i^{'},0}-(-1)^{\bar{i^{'}}}\frac{\lambda_{i^{'}}}{p^{a}}\xi_{j^{'}}^{j^{'},0}]=0$.
This gives us $a_{i,j,k}=0$ for all $1\leq i\neq j\leq r+s$ and
thus the element $x$ is reduced to the form 
\[
x=\sum_{\substack{0<k\leq\lambda_{i}-1\\
1\leq i\leq r+s
}
}b_{i,k}\xi_{i}^{i,k}+\sum_{1<i<r+s}c_{i,0}\left(\frac{\lambda_{i+1}}{p^{a}}\xi_{i}^{i,0}-(-1)^{\bar{i}}\frac{\lambda_{i}}{p^{a}}\xi_{i+1}^{i+1,0}\right).
\]

Next we consider an element $\xi_{i^{'}}^{j^{'},0}\in\mathfrak{g}^{e}$
for some $i^{'}<j^{'}$. Using a similar argument as above to $[x,\xi_{i^{'}}^{j^{'},0}]$
we obtain that $c_{i,0}=0$ for $1<i<r+s$ and $b_{j^{'},k}=b_{i^{'},k}$
for all $1\leq i^{'}\neq j^{'}\leq r+s$ and $0<k\leq\lambda_{1}-1$.
Note that $e^{k}=\sum_{i=1}^{r+s}\xi_{i}^{i,k}$. Therefore, we deduce
that $x\in\mathfrak{z}(\mathfrak{g}^{e})$ if and only if $x\in\langle e,\dots,e^{\lambda_{1}-1}\rangle$.

2. When $r+s=2$, i.e. $\lambda=(\lambda_{1}|\lambda_{2})$ and $\lambda_{1}=\lambda_{2}=n>1$.
Then an element $y\in\mathfrak{g}^{e}$ can be written as $y=\sum_{1\leq i\neq j\leq2,0\leq k\leq n-1}a_{i,j,k}\xi_{i}^{j,k}+\sum_{0<k\leq n-1,1\leq i\leq2}b_{i,k}\xi_{i}^{i,k}$
for $a_{i,j,k},b_{i,k}\in\mathbb{K}$. Now assume $y\in\mathfrak{z}(\mathfrak{g}^{e})$.
Applying a similar argument as in case $1$ to $[y,\xi_{1}^{1,1}]$
and $[y,\xi_{2}^{1,0}]$, we deduce that $a_{1,2,k}=a_{2,1,k}=0$
for $0\leq k\leq n-1$ and $b_{1,k}=b_{2,k}$ for $0<k\leq n-1$. 

In conclusion, we have that $\mathfrak{z}(\mathfrak{g}^{e})=\langle e,\dots,e^{\lambda_{1}-1}\rangle$.
\end{proof}

\section{Centralizer and Centre of centralizer of nilpotent elements in ortho-symplectic
Lie superalgebras\label{sec:osp}}

\subsection{Centralizer of nilpotent even elements in ortho-symplectic Lie superalgebras\label{subsec:centralizer-osp}}

Let $V=V_{\bar{0}}\oplus V_{\bar{1}}$ be a finite-dimensional vector
superspace over $\mathbb{K}$ such that $\dim V_{\bar{0}}=m$ and
$\dim V_{\bar{1}}=2n$. Let $B$ be a non-degenerate even supersymmetric
bilinear form on $V$, that is, $V_{\bar{0}}$ and $V_{\bar{1}}$
are orthogonal, the restriction of $B$ to $V_{\bar{0}}$ is symmetric,
and the restriction of $B$ to $V_{\bar{1}}$ is skew-symmetric. Then
the ortho-symplectic Lie superalgebra $\mathfrak{g}=\mathfrak{osp}(V)=\mathfrak{osp}(m|2n)$
is defined to be $\mathfrak{osp}(V)=\mathfrak{g}_{\bar{0}}\oplus\mathfrak{g}_{\bar{1}}$
where 
\[
\mathfrak{g}_{\bar{i}}=\{x\in\mathfrak{gl}(V)_{\bar{i}}:B(x(v),w)=-(-1)^{\bar{i}\bar{x}}B(v,x(w))\}
\]
 for homogeneous $v,w\in V$ and $\bar{i}\in\mathbb{Z}_{2}$.

Let $e\in\mathfrak{g}_{\bar{0}}$ be nilpotent with the corresponding
Jordan type $\lambda$ such that $\lambda=(\lambda_{1},\dots,\lambda_{r+s})$
is defined in a similar way to (\ref{eq:gl(m,n)-partition 1}) where
$\sum_{\bar{i}=\bar{0}}\lambda_{i}=m$ and $\sum_{\bar{i}=\bar{1}}\lambda_{i}=2n$.
In this subsection, we describe the structure of the centralizer $\mathfrak{g}^{e}$.
In case of zero characteristic, a basis for $\mathfrak{g}_{\mathbb{C}}^{e}$
was given in \cite[Subsection 4.5]{han-classical}, in particular,
it contains all elements in the set $S$ in (\ref{eq:basis-osp}).

According to \cite[Subsection 4.5]{han-classical}, there exists an
involution $i\mapsto i^{*}$ on the set $\{1,\dots,r+s\}$ such that
$i^{*}\in\{i-1,i,i+1\}$ for all $1\leq i\leq r+s$, $\lambda_{i}=\lambda_{i^{*}}$
and $\bar{i^{*}}=\bar{i}$. Note that $i=i^{*}$ if $\bar{i}=\bar{0}$
and $(-1)^{\lambda_{i}}=-1$ or if $\bar{i}=\bar{1}$ and $(-1)^{\lambda_{i}}=1$.
Recall that basis elements $\xi_{i}^{j,k}$ for $\mathfrak{gl}(m|2n)^{e}$
are defined in (\ref{eq:basis gl}). Let $S=S_{0}\cup S_{1}\cup S_{2}$
be a subset of $\mathfrak{g}$ given by
\begin{equation}
S_{0}=\{\xi_{i}^{i^{*},\lambda_{i}-1-k}:0\leq k\leq\lambda_{i}-1,\bar{i}=\bar{0},k\text{ is odd}\},\label{eq:basis-osp}
\end{equation}
\[
S_{1}=\{\xi_{i}^{i^{*},\lambda_{i}-1-k}:0\leq k\leq\lambda_{i}-1,\bar{i}=\bar{1},k\text{ is even}\},
\]
\[
\text{and }S_{2}=\{\xi_{i}^{j,\lambda_{j}-1-k}+\varepsilon_{i,j,k}\xi_{j^{*}}^{i^{*},\lambda_{i}-1-k}:0\leq k\leq\min\{\lambda_{i},\lambda_{j}\}-1,i^{*}\neq j,\varepsilon_{i,j,k}\in\{\pm1\}\}.
\]

\begin{doublespace}
\noindent We have that $\varepsilon_{i,j,k}=(-1)^{\lambda_{j}-k-\bar{x}\bar{i}}\theta_{j}\theta_{i}$
where $\bar{x}\in\{\bar{0},\bar{1}\}$ is the parity of $\xi_{i}^{j,\lambda_{j}-1-k}$
and $\theta_{j},\theta_{i}\in\{\pm1\}$ can be determined explicitly
by \cite[Section 3.2]{Janzten}. Note that an element $\xi_{i}^{j,\lambda_{j}-1-k}+\varepsilon_{i,j,k}\xi_{j^{*}}^{i^{*},\lambda_{i}-1-k}$
lies in $\mathfrak{g}_{\bar{0}}^{e}$ (resp. $\mathfrak{g}_{\bar{1}}^{e}$
) if $\bar{i}=\bar{j}$ (resp. $\bar{i}\neq\bar{j}$).
\end{doublespace}
\begin{thm}
We have that $S$ is a basis for $\mathfrak{g}^{e}$.
\end{thm}

\begin{proof}
We know that $e=\sum_{i=1}^{r+s}\xi_{i}^{i,1}$. Commutators between
$e$ and elements in $S$ over $\mathbb{Q}$ have been calculated
explicitly in \cite[Subsection 4.5]{han-classical}, they all equal
to zeros thus they must also be zeros in case of $\mathbb{K}$. Hence,
we deduce that $S$ is a subset of $\mathfrak{g}^{e}$. Based on the
conditions on indices $i,j,k$ and that $\xi_{i}^{j,\lambda_{j}-1-k}$
are linearly independent in $\mathfrak{gl}(m|2n)$, for $c_{i,i^{*},k},c_{i,j,k}\in\mathbb{K}$,
$\sum c_{i,i^{*},k}\xi_{i}^{i^{*},\lambda_{i}-1-k}+\sum_{i^{*}\neq j}c_{i,j,k}(\xi_{i}^{j,\lambda_{j}-1-k}+\varepsilon_{i,j,k}\xi_{j^{*}}^{i^{*},\lambda_{i}-1-k})=0$
implies that $c_{i,i^{*},k}=c_{i,j,k}=0$. Thus we have that elements
$\xi_{i}^{i^{*},\lambda_{i}-1-k}$ and $\xi_{i}^{j,\lambda_{j}-1-k}+\varepsilon_{i,j,k}\xi_{j^{*}}^{i^{*},\lambda_{i}-1-k}$
in $S$ are linearly independent. By \cite[Proposition 20]{han-classical},
we know that $|S|=\dim\mathfrak{g}_{\mathbb{C}}(0)+\dim\mathfrak{g}_{\mathbb{C}}(1)$.
According to Lemma \ref{lem:g(j)}, we have that $\dim\mathfrak{g}_{\mathbb{C}}(j)=\dim\mathfrak{g}(j)$
for all $j\in\mathbb{Z}$ and thus $|S|=\dim\mathfrak{g}(0)+\dim\mathfrak{g}(1)$.
Furthermore, recall that $\dim\mathfrak{g}^{e}=\dim\mathfrak{g}(0)+\dim\mathfrak{g}(1)$
by Theorem \ref{thm:g^e=00003Ddimg(0)+dimg(1)}. Therefore, we deduce
that $|S|=\dim\mathfrak{g}^{e}$ and $S$ is a basis for $\mathfrak{g}^{e}$. 
\end{proof}

\subsection{Centre of centralizer of nilpotent even elements in ortho-symplectic
Lie superalgebras}

\begin{singlespace}
\noindent A basis of $\mathfrak{z}(\mathfrak{g}_{\mathbb{C}}^{e})$
is given in \cite[Section 4]{han-classical} and we apply a similar
argument to $\mathfrak{g}$ in this subsection. In particular, we
prove Theorem \ref{thm:z(g^e)-osp} which gives a basis of $\mathfrak{z}(\mathfrak{g}^{e})$,
though note that some explicit calculations are omitted as they are
the same as in \cite[Section 4]{han-classical}.
\end{singlespace}

Before we start the proof, we need the following notation. In this
subsection, we use an alternative notation for $\lambda$ and rewrite
\begin{equation}
\lambda=(\lambda_{1},\dots,\lambda_{a},\lambda_{a+1},\lambda_{-(a+1)},\dots,\lambda_{b},\lambda_{-b})\label{eq:jordan type group}
\end{equation}
 where $\lambda_{1},\dots,\lambda_{a}$ are the parts with odd multiplicity,
$\lambda_{1}>\lambda_{2}>\dots>\lambda_{a}$ and $\lambda_{a+1}=\lambda_{-(a+1)}\geq\dots\geq\lambda_{b}=\lambda_{-b}$.
We define $\bar{i}\in\{\bar{0},\bar{1}\}$ such that for $c\in\mathbb{Z}$,
we have $\ensuremath{\left|\{i:\lambda_{i}=c,\bar{i}=\bar{0}\}\right|}=\ensuremath{\left|\{j:p_{j}=c\}\right|}$
and $\left|\{i:\lambda_{i}=c,\bar{i}=\bar{1}\}\right|=\ensuremath{\left|\{j:q_{j}=c\}\right|}$
for some $j$. Then there exists $\{v_{1},\dots,v_{a},v_{a+1},v_{-(a+1)},\dots,v_{b},v_{-b}\}\in V$
such that our vector superspace $V$ can be decomposed as $V=V_{1}\oplus V_{2}$
where $V_{1}=\langle e^{j}v_{i}:1\leq i\leq a\rangle$ and $V_{2}=\langle e^{j}v_{i},e^{j}v_{-i}:a+1\leq i\leq b\rangle$.
Let $\mathfrak{g}_{1}=\mathfrak{osp}(V_{1})$ and $\mathfrak{g}_{2}=\mathfrak{osp}(V_{2})$,
we define $\bar{\mathfrak{g}}=\mathfrak{g}_{1}\oplus\mathfrak{g}_{2}$.
Then a nilpotent element $e\in\mathfrak{g}_{\bar{0}}$ can be viewed
as $e=e_{1}+e_{2}$ where $e_{i}\in\mathfrak{g}_{i}$ and the Jordan
type of $e_{1}$ (resp. $e_{2}$) is $(\lambda_{1},\dots,\lambda_{a})$
(resp. $(\lambda_{a+1},\lambda_{-(a+1)},\dots,\lambda_{b},\lambda_{-b})$).
We have that $\mathfrak{z}(\mathfrak{g}^{e})\subseteq\mathfrak{z}(\mathfrak{g}_{1}^{e_{1}})\oplus\mathfrak{z}(\mathfrak{g}_{2}^{e_{2}})$
by \cite[Subsection 4.8]{han-classical}. 
\begin{thm}
\label{thm:z(g^e)-osp}Let $A=\{e^{l}:l\text{ is odd and }0\leq l\leq\min\{\lambda_{1},\lambda_{a+1}\}-1\}$.
Then $A$ is a basis for $\mathfrak{z}(\mathfrak{g}^{e})$ except
for the following two cases:

$\bullet$ If $a\geq3$, $\lambda_{2}>\lambda_{a+1}$ and $\bar{i}=\bar{0}$
for $i=1,2$ or $a=2,\bar{i}\neq\bar{j}$ for $(i,j)=(1,2)$, then
we have that $A\cup\{\xi_{1}^{2,\lambda_{2}-1}-\xi_{2}^{1,\lambda_{1}-1}\}$
is a basis for $\mathfrak{z}(\mathfrak{g}^{e})$;

$\bullet$ If $\lambda_{1}<\lambda_{a+1}$, $\lambda_{a+1}>\lambda_{a+2}$,
$\bar{i}=\bar{0}$ for $i=a+1$ and $\lambda_{a+1}$ is odd, we have
that $A\cup\{\xi_{a+1}^{a+1,\lambda_{a+1}-1}-\xi_{-(a+1)}^{-(a+1),\lambda_{a+1}-1}\}$
is a basis for $\mathfrak{z}(\mathfrak{g}^{e})$.
\end{thm}

\begin{proof}
\begin{singlespace}
\noindent Step 1: Determine $\mathfrak{z}(\mathfrak{g}_{1}^{e_{1}})$.
\end{singlespace}

According to Subsection \ref{subsec:centralizer-osp}, we know that
an element $x_{1}\in\mathfrak{z}(\mathfrak{g}_{1}^{e_{1}})$ can be
written as 
\[
x_{1}=\sum_{\substack{0\leq k\leq\lambda_{i}-1,\lambda_{i}-k\text{ is even}\\
1\leq i\leq a
}
}c_{i,k}\xi_{i}^{i,\lambda_{i}-1-k}+\sum_{\substack{0\leq k\leq\lambda_{j}-1\\
1\leq i<j\leq a
}
}c_{i,j,k}(\xi_{i}^{j,\lambda_{j}-1-k}+\varepsilon_{i,j,k}\xi_{j}^{i,\lambda_{i}-1-k})
\]
 where $c_{i,k},c_{i,j,k}\in\mathbb{K}$. 

Suppose $a\geq3$. By considering $[\xi_{t}^{t,1},x_{1}]=0$ for $1\leq t\leq a$,
we obtain that $c_{i,j,k}=0$ for $1\leq k\leq\lambda_{j}-1$. For
$1\leq t<h\leq a$, the commutator $[\xi_{t}^{h,0}+\varepsilon_{t,h,\lambda_{h}-1}\xi_{h}^{t,\lambda_{t}-\lambda_{h}},x_{1}]=0$
implies that $c_{t,k}=c_{h,k}$ for all $1\leq t<h\leq a$ and $c_{i,j,0}=0$
for all $1\leq i<j\leq a$ except when $\bar{i}=\bar{0}$ for $i=1,2$.
Hence, we deduce that $\mathfrak{z}(\mathfrak{g}_{1}^{e_{1}})=\langle\sum_{i=1}^{a}\xi_{i}^{i,l}:l\text{ is odd and }1\leq l\leq\lambda_{1}-1\rangle$
except when $\bar{i}=\bar{0}$ for $i=1,2$, in which case we have
that $\xi_{1}^{2,\lambda_{2}-1}-\xi_{2}^{1,\lambda_{1}-1}$ commutes
with all basis elements in $\mathfrak{g}^{e}$. 

Applying a similar argument as above to the case in which $a=2$,
we obtain that $c_{1,2,k}=0$ for $1\leq k\leq2n-1$, $c_{1,k}=c_{2,k}$
for all $k$ and $\xi_{1}^{2,\lambda_{2}-1}-\xi_{2}^{1,\lambda_{1}-1}$
commutes with all basis elements in $\mathfrak{g}^{e}$. 

Therefore, we have that $\mathfrak{z}(\mathfrak{g}_{1}^{e_{1}})=\langle e_{1}^{l}:l\text{ is odd and }1\leq l\leq\lambda_{1}-1\rangle$
except when $a\geq3$ and $\bar{i}=\bar{0}$ for $i=1,2$ or $a=2,\bar{i}\neq\bar{j}$
for $(i,j)=(1,2)$, in which case we have that $\mathfrak{z}(\mathfrak{g}_{1}^{e_{1}})=\langle e_{1}^{l}:l\text{ is odd and }1\leq l\leq\lambda_{1}-1\rangle\oplus\langle\xi_{1}^{2,\lambda_{2}-1}-\xi_{2}^{1,\lambda_{1}-1}\rangle$.

Step 2: Determine $\mathfrak{z}(\mathfrak{g}_{2}^{e_{2}})$.

Based on arguments in \cite[Subsection 4.6]{han-classical}, we have
that $\mathfrak{z}(\mathfrak{g}_{2}^{e_{2}})\subseteq(\mathfrak{g}_{2}^{e_{2}})^{\mathfrak{h}_{2}^{e_{2}}}$
where $\mathfrak{h}_{2}$ is a Cartan subalgebra of $\mathfrak{g}_{2}$
and $\mathfrak{h}_{2}^{e_{2}}=\langle h_{i}=\xi_{i}^{i,0}-\xi_{-i}^{-i,0}:a+1\leq i\leq b\rangle$.
By computing commutators between an element $h\in\mathfrak{h}_{2}^{e_{2}}$
and basis elements of $\mathfrak{g}_{2}^{e_{2}}$, we deduce that
$(\mathfrak{g}_{2}^{e_{2}})^{\mathfrak{h}_{2}^{e_{2}}}=\langle\xi_{j}^{j,\lambda_{j}-1-k}+\varepsilon_{j,j,k}\xi_{-j}^{-j,\lambda_{j}-1-k}:a+1\leq j\leq b\rangle$.
Now let $x_{2}\in\mathfrak{z}(\mathfrak{g}_{2}^{e_{2}})$ such that
$x_{2}=\sum_{j,k}c_{j,k}(\xi_{j}^{j,\lambda_{j}-1-k}+\varepsilon_{j,j,k}\xi_{-j}^{-j,\lambda_{j}-1-k})$
where $c_{j,k}\in\mathbb{K}$. Note that for $a+1\leq i\leq b$, an
element $\xi_{i}^{-i,0}\in\mathfrak{g}_{2}^{e_{2}}$ (resp. $\xi_{i}^{-i,1}\in\mathfrak{g}_{2}^{e_{2}}$)
if $\lambda_{i}$ is even and $\bar{i}=\bar{0}$ (resp. $\lambda_{i}$
is odd and $\bar{i}=\bar{1}$), then $[\xi_{i}^{-i,0},x_{2}]=0$ (resp.
$[\xi_{i}^{-i,1},x_{2}]=0$) implies that $c_{j,k}=0$ if $\lambda_{j}-k$
is odd and $k\neq0$. For $\lambda_{j}-k$ is odd and $k=0$, if there
exists some $\lambda_{i}$ such that $a+1\leq i\leq b$ and $\lambda_{i}\geq\lambda_{j}$,
by considering $[\xi_{i}^{j,0}+\varepsilon_{i,j,\lambda_{j}-1}\xi_{-j}^{-i,\lambda_{i}-\lambda_{j}},x_{2}]=0$
we have that $c_{j,0}=0$. Hence, we deduce that $\mathfrak{z}(\mathfrak{g}_{2}^{e_{2}})\subseteq\langle\xi_{j}^{j,\lambda_{j}-1-k}+\xi_{-j}^{-j,\lambda_{j}-1-k}:a+1\leq j\leq b,\lambda_{j}-k\text{ is even}\rangle$
except when $\lambda_{a+1}$ is odd, $\lambda_{a+1}>\lambda_{i}$
for all $a+2\leq|i|\leq b$ and $\overline{a+1}=\bar{0}$, in which
case $\xi_{a+1}^{a+1,\lambda_{a+1}-1}-\xi_{-(a+1)}^{-(a+1),\lambda_{a+1}-1}$
commutes with all basis elements in $\mathfrak{g}_{2}^{e_{2}}$. Furthermore,
for $a+1\leq i<t\leq b$, the commutator $[\xi_{i}^{t,0}+\varepsilon_{i,t,\lambda_{t}-1}\xi_{-t}^{-i,\lambda_{i}-\lambda_{t}},x_{2}]=0$
gives us $c_{i,k}=c_{t,k}$ whenever all $\lambda_{i}-k$ and $\lambda_{t}-k$
are even.

Therefore, we have that $\mathfrak{z}(\mathfrak{g}_{2}^{e_{2}})=\langle e_{2}^{l}:l\text{ is odd and }1\leq l\leq\lambda_{a+1}-1\rangle$
except when $\lambda_{a+1}$ is odd, $\lambda_{a+1}>\lambda_{i}$
for all $a+2\leq|i|\leq b$ and $\overline{a+1}=\bar{0}$, in which
case we have that $\mathfrak{z}(\mathfrak{g}_{2}^{e_{2}})=\langle e_{2}^{l}:l\text{ is odd and }1\leq l\leq\lambda_{a+1}-1\rangle\oplus\langle\xi_{a+1}^{a+1,\lambda_{a+1}-1}-\xi_{-(a+1)}^{-(a+1),\lambda_{a+1}-1}\rangle$.

Step 3: Determine $\mathfrak{z}(\mathfrak{g}^{e})$.

We first consider the case when $\mathfrak{z}(\mathfrak{g}_{1}^{e_{1}})=\langle e_{1}^{l}:l\text{ is odd and }1\leq l\leq\lambda_{1}-1\rangle$
and $\mathfrak{z}(\mathfrak{g}_{2}^{e_{2}})=\langle e_{2}^{l}:l\text{ is odd and }1\leq l\leq\lambda_{a+1}-1\rangle$.
An element $x\in\mathfrak{z}(\mathfrak{g}^{e})$ is of the form 
\[
x=\sum_{l\text{ is odd;}l=1}^{\lambda_{1}-1}a_{l}\left(\sum_{t=1}^{a}\xi_{t}^{t,l}\right)+\sum_{l\text{ is odd;}k=1}^{\lambda_{a+1}-1}b_{l}\left(\sum_{t=a+1}^{b}(\xi_{t}^{t,l}+\xi_{-t}^{-t,l})\right)
\]
 for $a_{l},b_{l}\in\mathbb{K}$. If $\lambda_{1}\geq\lambda_{a+1}$
(resp. $\lambda_{1}<\lambda_{a+1}$), by letting $[\xi_{1}^{a+1,0}+\varepsilon_{1,a+1,\lambda_{a+1}-1}\xi_{-(a+1)}^{1,\lambda_{1}-\lambda_{a+1}},x]=0$
(resp. $[\xi_{1}^{a+1,\lambda_{a+1}-\lambda_{1}}+\varepsilon_{1,a+1,\lambda_{1}-1}\xi_{-(a+1)}^{1,0},x]=0$)
we obtain that $a_{l}=b_{l}$ for all odd $l$. Hence, in this case
we deduce that $A$ is a basis for $\mathfrak{z}(\mathfrak{g}^{e})$. 

It remains to consider special cases. (i) When $\xi_{1}^{2,\lambda_{2}-1}-\xi_{2}^{1,\lambda_{1}-1}\in\mathfrak{z}(\mathfrak{g}_{1}^{e_{1}})$
and $\xi_{a+1}^{a+1,\lambda_{a+1}-1}-\xi_{-(a+1)}^{-(a+1),\lambda_{a+1}-1}\notin\mathfrak{z}(\mathfrak{g}_{2}^{e_{2}})$,
an element $y\in\mathfrak{z}(\mathfrak{g}^{e})$ is of the form $y=x+c_{1,2,0}(\xi_{1}^{2,\lambda_{2}-1}-\xi_{2}^{1,\lambda_{1}-1})$.
Applying a similar argument as above to $y$ gives us $a_{l}=b_{l}$
for all odd $l$ and $c_{1,2,0}=0$ except for $\lambda_{2}>\lambda_{a+1}$,
in which case $\xi_{1}^{2,\lambda_{2}-1}-\xi_{2}^{1,\lambda_{1}-1}$
commutes with all basis elements of $\mathfrak{g}^{e}$. 

(ii) When $\xi_{1}^{2,\lambda_{2}-1}-\xi_{2}^{1,\lambda_{1}-1}\notin\mathfrak{z}(\mathfrak{g}_{1}^{e_{1}})$
and $\xi_{a+1}^{a+1,\lambda_{a+1}-1}-\xi_{-(a+1)}^{-(a+1),\lambda_{a+1}-1}\in\mathfrak{z}(\mathfrak{g}_{2}^{e_{2}})$,
an element $z\in\mathfrak{z}(\mathfrak{g}^{e})$ is of the form $z=x+c_{a+1,a+1,0}(\xi_{a+1}^{a+1,\lambda_{a+1}-1}-\xi_{-(a+1)}^{-(a+1),\lambda_{a+1}-1})$.
Applying a similar argument as above to $z$ gives us $a_{l}=b_{l}$
for all odd $l$ and $c_{a+1,a+1,0}=0$ except for $\lambda_{1}<\lambda_{a+1}$,
in which case $\xi_{a+1}^{a+1,\lambda_{a+1}-1}-\xi_{-(a+1)}^{-(a+1),\lambda_{a+1}-1}$
commutes with all basis elements of $\mathfrak{g}^{e}$.

(iii) When $\xi_{1}^{2,\lambda_{2}-1}-\xi_{2}^{1,\lambda_{1}-1}\in\mathfrak{z}(\mathfrak{g}_{1}^{e_{1}})$
and $\xi_{a+1}^{a+1,\lambda_{a+1}-1}-\xi_{-(a+1)}^{-(a+1),\lambda_{a+1}-1}\in\mathfrak{z}(\mathfrak{g}_{2}^{e_{2}})$,
applying a similar argument as above we also have that $\mathfrak{z}(\mathfrak{g}^{e})=\langle A\rangle\oplus\langle\xi_{1}^{2,\lambda_{2}-1}-\xi_{2}^{1,\lambda_{1}-1}\rangle$
for $\lambda_{2}>\lambda_{a+1}$ and $\mathfrak{z}(\mathfrak{g}^{e})=\langle A\rangle\oplus\langle\xi_{a+1}^{a+1,\lambda_{a+1}-1}-\xi_{-(a+1)}^{-(a+1),\lambda_{a+1}-1}\rangle$
for $\lambda_{1}<\lambda_{a+1}$.
\end{proof}

\section{Centralizer and centre of centralizer of nilpotent elements in exceptional
Lie superalgebras\label{sec:exceptional}}

\noindent In this section, we fix the following notation. For an element
$x=x_{\bar{0}}+x_{\bar{1}}\in\mathfrak{g}$, we have that $[e,x]=[e,x_{\bar{0}}]+[e,x_{\bar{1}}]$
where $[e,x_{\bar{0}}]\in\mathfrak{g}_{\bar{0}}$ and $[e,x_{\bar{1}}]\in\mathfrak{g}_{\bar{1}}$.
This implies that $\mathfrak{g}^{e}=\mathfrak{g}_{\bar{0}}^{e}\oplus\mathfrak{g}_{\bar{1}}^{e}$.
Similarly, let us denote $\mathfrak{z}(\mathfrak{g}^{e})=\mathfrak{z}_{\bar{0}}\oplus\mathfrak{z}_{\bar{1}}$
and $\mathfrak{z}_{\bar{i}}(j)=\mathfrak{z}_{\bar{i}}\cap\mathfrak{g}(j)$.
We often work with $\mathfrak{sl}_{2}(\mathbb{K})$ so we fix notation
$E,H,F$ for its basis elements where 
\[
E=\begin{pmatrix}0 & 1\\
0 & 0
\end{pmatrix},H=\begin{pmatrix}1 & 0\\
0 & -1
\end{pmatrix},F=\begin{pmatrix}0 & 0\\
1 & 0
\end{pmatrix}.
\]

\subsection{Centralizer of even nilpotent elements in $D(2,1;\alpha)$ and its
centre\label{subsec:D(2,1;)}}

The Lie superalgebras $D(2,1;\alpha)$ with $\alpha\in\mathbb{K}\backslash\{0,-1\}$
form a one-parameter family of seventeen-dimensional Lie superalgebras.
Note that $D(2,1;\alpha)$ is also denoted by $\Gamma(\sigma_{1},\sigma_{2},\sigma_{3})$
in \cite{M.Scheunert1976} where $\sigma_{1},\sigma_{2},\sigma_{3}\in\mathbb{K}\backslash\{0,1\}$
and $\sigma_{1}+\sigma_{2}+\sigma_{3}=0$. According to \cite[Lemma 5.5.16]{Musson2012},
if there is another triple $(\sigma_{1}^{'},\sigma_{2}^{'},\sigma_{3}^{'})$
such that $\sigma_{1}^{'}+\sigma_{2}^{'}+\sigma_{3}^{'}=0$ and $\Gamma(\sigma_{1}^{'},\sigma_{2}^{'},\sigma_{3}^{'})\cong\Gamma(\sigma_{1},\sigma_{2},\sigma_{3})$,
then there exist a permutation $\rho$ of $\{1,2,3\}$ and a nonzero
complex number $c$ such that $\sigma_{i}^{'}$ can be obtained by
$\sigma_{i}^{'}=c\sigma_{\rho(i)}$. For any $\alpha\in\mathbb{K}\backslash\{0,-1\}$,
we have $D(2,1;\alpha)=\Gamma(1+\alpha,-1,-\alpha)\cong\Gamma(-\alpha,-1,1+\alpha)\cong\Gamma(\frac{1+\alpha}{\alpha},-1,-\frac{1}{\alpha})$.

Let $\mathfrak{g}=\mathfrak{g}_{\bar{0}}\oplus\mathfrak{g}_{\bar{1}}=D(2,1;\alpha)$
with $\alpha\in\mathbb{K}\backslash\{0,-1\}$. By definition the even
part is $\mathfrak{g}_{\bar{0}}=\mathfrak{sl}_{2}(\mathbb{K})\oplus\mathfrak{sl}_{2}(\mathbb{K})\oplus\mathfrak{sl}_{2}(\mathbb{K})$
and the odd part is $\mathfrak{g}_{\bar{1}}=V_{1}\otimes V_{2}\otimes V_{3}$
where $V_{i}$, $i=1,2,3$ is the standard two-dimensional $\mathfrak{sl}_{2}(\mathbb{K})$-module
for the $i$th summand in $\mathfrak{g}_{\bar{0}}$ with basis elements
$v_{1}^{i}=(1,0)^{t}$ and $v_{-1}^{i}=(0,1)^{t}$. Note that the
map $\mathfrak{g}_{\bar{0}}\times\mathfrak{g}_{\bar{0}}\rightarrow\mathfrak{g}_{\bar{0}}$
is a Lie bracket and the Lie superbracket $[\cdotp,\cdotp]:\mathfrak{g}_{\bar{0}}\times\mathfrak{g}_{\bar{1}}\rightarrow\mathfrak{g}_{\bar{1}}$
is defined by $[x_{1}\oplus x_{2}\oplus x_{3},v_{i}^{1}\otimes v_{j}^{2}\otimes v_{k}^{3}]=x_{1}v_{i}^{1}\otimes v_{j}^{2}\otimes v_{k}^{3}+v_{i}^{1}\otimes x_{2}v_{j}^{2}\otimes v_{k}^{3}+v_{i}^{1}\otimes v_{j}^{2}\otimes x_{3}v_{k}^{3}$
for $x_{1}\oplus x_{2}\oplus x_{3}\in\mathfrak{g}_{\bar{0}}$, $v_{i}^{1}\otimes v_{j}^{2}\otimes v_{k}^{3}\in\mathfrak{g}_{\bar{1}}$.
The Lie superbracket $[\cdotp,\cdotp]:\mathfrak{g}_{\bar{1}}\times\mathfrak{g}_{\bar{1}}\rightarrow\mathfrak{g}_{\bar{0}}$
is given for example in \cite[Section 4.2]{han-exp} and it depends
on $\sigma_{i}$ for $i=1,2,3$. 

In this subsection, we denote $E_{1}=(E,0,0)$, $E_{2}=(0,E,0)$ and
$E_{3}=(0,0,E)$. Similarly, let $F_{1}=(F,0,0)$, $F_{2}=(0,F,0)$,
$F_{3}=(0,0,F)$, $H_{1}=(H,0,0)$, $H_{2}=(0,H,0)$ and $H_{3}=(0,0,H)$.
We write $v_{i,j,k}$ for $v_{i}^{1}\otimes v_{j}^{2}\otimes v_{k}^{3}$
where $i,j,k\in\{\pm1\}$. Then a basis for $\mathfrak{g}_{\bar{0}}$
is $\{E_{i},H_{i},F_{i}:i=1,2,3\}$ and a basis for $\mathfrak{g}_{\bar{1}}$
is $\{v_{i,j,k}:i,j,k=\pm1\}$.

\begin{singlespace}
Since representatives of nilpotent orbits in $\mathfrak{sl}_{2}(\mathbb{K})$
up to conjugation by $\mathrm{SL}_{2}(\mathbb{K})$ are $0$ and $E$,
we have that representatives of nilpotent orbits $e\in\mathfrak{g}_{\bar{0}}$
are $e=0,E_{1},E_{2},E_{3},E_{1}+E_{2},E_{2}+E_{3},E_{1}+E_{3},E_{1}+E_{2}+E_{3}$.
We give the cocharacter $\tau$ associated to $e$ and basis elements
for $\mathfrak{g}^{e}$ and $\mathfrak{z}\left(\mathfrak{g}^{e}\right)$
for $e=0,E_{1},E_{1}+E_{2},E_{1}+E_{2}+E_{3}$ in Table \ref{tab:D(2,1)}.
Note that in the second column of Table \ref{tab:D(2,1)}, these $\alpha_{i}=2\beta_{i}$,
$i=1,2,3$ are even roots in the root system for $\mathfrak{g}$ as
given in Subsection \ref{sec:highest root}. After Table \ref{tab:D(2,1)},
we explain our explicit calculation for the case $e=E_{1}+E_{2}+E_{3}$.
The cases when $e=0$, $E_{1}$, $E_{1}+E_{2}$ are obtained using
a similar approach. Note that analysis for cases $e=E_{2}$, $e=E_{3}$
are similar to $e=E_{1}$ and analysis for cases $e=E_{1}+E_{3}$,
$e=E_{2}+E_{3}$ are similar to $e=E_{1}+E_{2}$, for which are omitted
in this paper.
\end{singlespace}

\begin{table}[H]
\noindent \begin{centering}
\begin{tabular}{|>{\centering}p{2.5cm}|>{\centering}p{2.8cm}|>{\centering}p{5cm}|>{\centering}p{1.2cm}|>{\centering}p{1.3cm}|}
\hline 
Representatives of nilpotent orbits $e\in\mathfrak{g}_{\bar{0}}$ & Cocharacter $\tau$ & $\mathfrak{g}^{e}=\mathfrak{g}_{\bar{0}}^{e}\oplus\mathfrak{g}_{\bar{1}}^{e}$ & $\dim\left(\mathfrak{g}^{e}\right)$ & $\mathfrak{z}\left(\mathfrak{g}^{e}\right)$\tabularnewline
\hline 
\hline 
$0$ & $0$ & $\mathfrak{g}$ & $17$ & $\langle e\rangle=0$\tabularnewline
\hline 
$E_{1}$ & $h_{\alpha_{1}}(t)$ & $\langle E_{1},E_{2},H_{2},F_{2},E_{3},H_{3},F_{3}\rangle\oplus\langle v_{1,j,k}:j,k=\pm1\rangle$ & $11$ & $\langle e\rangle$\tabularnewline
\hline 
$E_{1}+E_{2}$ & $h_{\alpha_{1}}(t)h_{\alpha_{2}}(t)$ & $\langle E_{1},E_{2},E_{3},H_{3},F_{3}\rangle\oplus\langle v_{1,1,1},v_{1,1,-1},v_{1,-1,1}-v_{-1,1,1},v_{1,-1,-1}-v_{-1,1,-1}\rangle$ & $9$ & $\langle e\rangle$\tabularnewline
\hline 
$E_{1}+E_{2}+E_{3}$ & $h_{\alpha_{1}}(t)h_{\alpha_{2}}(t)h_{\alpha_{3}}(t)$ & $\langle E_{1},E_{2},E_{3}\rangle\oplus\langle v_{1,1,1},v_{1,1,-1}-v_{-1,1,1},v_{1,-1,1}-v_{-1,1,1}\rangle$ & $6$ & $\langle e,v_{1,1,1}\rangle$\tabularnewline
\hline 
\end{tabular}
\par\end{centering}
\caption{\label{tab:D(2,1)}$\tau$, $\mathfrak{g}^{e}$ and $\mathfrak{z}\left(\mathfrak{g}^{e}\right)$}

\end{table}

Now we give our calculation explicitly for the case $e=E_{1}+E_{2}+E_{3}$.
Since $\mathfrak{g}_{\bar{0}}=\mathfrak{sl}_{2}(\mathbb{K})\oplus\mathfrak{sl}_{2}(\mathbb{K})\oplus\mathfrak{sl}_{2}(\mathbb{K})$,
we easily compute that $\mathfrak{g}_{\bar{0}}^{e}=\langle E_{1},E_{2},E_{3}\rangle$.
To determine $\mathfrak{g}_{\bar{1}}^{e}$, we give the $\tau$-grading
on basis elements of $\mathfrak{g}_{\bar{1}}$ in the following table. 

\begin{table}[H]
\noindent \begin{centering}
\begin{tabular}{|c|c|}
\hline 
$\tau$-grading & basis elements of $\mathfrak{g}_{\bar{1}}$\tabularnewline
\hline 
\hline 
$3$ & $v_{1,1,1}$\tabularnewline
\hline 
$1$ & $v_{1,1,-1},v_{1,-1,1},v_{-1,1,1}$\tabularnewline
\hline 
$-1$ & $v_{1,-1,-1},v_{-1,1,-1},v_{-1,-1,1}$\tabularnewline
\hline 
$-3$ & $v_{-1,-1,-1}$\tabularnewline
\hline 
\end{tabular}
\par\end{centering}
\caption{The $\tau$-grading on basis elements of $\mathfrak{g}_{\bar{1}}$
when $e=E_{1}+E_{2}+E_{3}$}

\end{table}
 From the above table we have that $\mathfrak{g}_{\bar{1}}^{e}=\mathfrak{g}_{\bar{1}}^{e}(3)\oplus\mathfrak{g}_{\bar{1}}^{e}(1)\oplus\mathfrak{g}_{\bar{1}}^{e}(-1)\oplus\mathfrak{g}_{\bar{1}}^{e}(-3)$.
We know that $[e,\mathfrak{g}_{\bar{1}}(3)]\subseteq\mathfrak{g}_{\bar{1}}(5)$
and since $\mathfrak{g}_{\bar{1}}(5)=0$, we have that $\mathfrak{g}_{\bar{1}}^{e}(3)=\langle v_{1,1,1}\rangle$.
In order to determine $\mathfrak{g}_{\bar{1}}^{e}(1)$, we consider
$x=a_{1,1,-1}v_{1,1,-1}+a_{1,-1,1}v_{1,-1,1}+a_{-1,1,1}v_{-1,1,1}\in\mathfrak{g}_{\bar{1}}^{e}(1)$
for $a_{1,1,-1}$, $a_{1,-1,1}$, $a_{-1,1,1}\in\mathbb{K}$. Then
we compute $[e,x]=(a_{1,1,-1}+a_{1,-1,1}+a_{-1,1,1})v_{1,1,1}$. This
equals to zero implies that $a_{1,1,-1}+a_{1,-1,1}+a_{-1,1,1}=0$.
Thus we have that $\mathfrak{g}_{\bar{1}}^{e}(1)=\langle v_{1,1,1},v_{1,1,-1}-v_{-1,1,1},v_{1,-1,1}-v_{-1,1,1}\rangle$.
To determine $\mathfrak{g}_{\bar{1}}^{e}(-1)$, we consider $y=a_{1,-1,-1}v_{1,-1,-1}+a_{-1,1,-1}v_{-1,1,-1}+a_{-1,-1,1}v_{-1,-1,1}\in\mathfrak{g}_{\bar{1}}^{e}(-1)$
for $a_{1,-1,-1}$, $a_{-1,1,-1}$, $a_{-1,-1,1}\in\mathbb{K}$. We
have that $[e,y]=(a_{-1,1,-1}+a_{1,-1,-1})v_{1,1,-1}+(a_{-1,-1,1}+a_{1,-1,-1})v_{1,-1,1}+(a_{-1,-1,1}+a_{-1,1,-1})v_{-1,1,1},$
this equals to zero implies that $a_{1,-1,-1}=a_{-1,1,-1}=a_{-1,-1,1}=0$.
Hence, we deduce that $\mathfrak{g}_{\bar{1}}^{e}(-1)=0$. We also
obtain that $\mathfrak{g}_{\bar{1}}^{e}(-3)=0$ since $[e,v_{-1,-1,-1}]=v_{1,-1,-1}+v_{-1,1,-1}+v_{-1,-1,1}\neq0$.
Therefore, we obtain that $\mathfrak{g}^{e}=\langle E_{1},E_{2},E_{3}\rangle\oplus\langle v_{1,1,1},v_{1,1,-1}-v_{-1,1,1},v_{1,-1,1}-v_{-1,1,1}\rangle$.

To determine $\mathfrak{z}(\mathfrak{g}^{e})$, we let $x=b_{1}E_{1}+b_{2}E_{2}+b_{3}E_{3}\in\mathfrak{z}_{\bar{0}}$.
Then $[x,v_{1,1,-1}-v_{-1,1,1}]=0$ and $[x,v_{1,-1,1}-v_{-1,1,1}]=0$
provide us with $b_{1}=b_{2}=b_{3}$. Hence we obtain that $\mathfrak{z}_{\bar{0}}=\langle E_{1}+E_{2}+E_{3}\rangle=\langle e\rangle$.
Next we look at $\mathfrak{z}_{\bar{1}}=\mathfrak{z}_{\bar{1}}(3)\oplus\mathfrak{z}_{\bar{1}}(1)$.
Let $z=b_{4}(v_{1,1,-1}-v_{-1,1,1})+b_{5}(v_{1,-1,1}-v_{-1,1,1})\in\mathfrak{z}_{\bar{1}}(3)$
for $b_{4},b_{5}\in\mathbb{K}$. We calculate $[E_{2},z]=b_{5}v_{1,1,1}$
and $[E_{3},z]=b_{4}v_{1,1,1}$, these commutators are equal to zero
if and only if $b_{4}=b_{5}=0$. Hence, we deduce that $\mathfrak{z}_{\bar{1}}(1)=0$.
Now we consider $\mathfrak{z}_{\bar{1}}(3)$. Observe that $\mathfrak{g}^{e}(3+j)=0$
for all $j$ such that $\mathfrak{g}^{e}(j)\neq0$. Hence, we have
that $[w,v_{1,1,1}]=0$ for all $w\in\mathfrak{g}^{e}(j)\neq0$ and
thus $\mathfrak{z}_{\bar{1}}(3)=\langle v_{1,1,1}\rangle$. Therefore,
we obtain that $\mathfrak{z}(\mathfrak{g}^{e})=\langle e,v_{1,1,1}\rangle$.

\subsection{Centralizer of even nilpotent elements in $G(3)$ and its centre\label{subsec:G(3)}}

Let $\mathfrak{g}=\mathfrak{g}_{\bar{0}}\oplus\mathfrak{g}_{\bar{1}}=G(3)$.
The explicit description of the Lie superalgebra $G(3)$ is given
for example in \cite[Section 5]{han-exp}. Recall that the even part
$\mathfrak{g}_{\bar{0}}$ is a direct sum of Lie algebras $\mathfrak{sl}_{2}(\mathbb{K})$
and $G_{2}$. The odd part $\mathfrak{g}_{\bar{1}}=V_{2}\otimes V_{7}$
where $V_{2}=\langle v_{1}=(1,0)^{t},v_{-1}=(0,1)^{t}\rangle$ is
the standard two-dimensional representation for $\mathfrak{sl}_{2}(\mathbb{K})$
and $V_{7}$ is a seven-dimensional simple representation of $G_{2}$
with a basis $\{e_{3},e_{2},e_{1},e_{0},e_{-1},e_{-2},e_{-3}\}$.
Thus $\mathfrak{g}_{\bar{1}}$ has a basis $\{v_{i}\otimes e_{j}:i=\pm1,j=0,\pm1,\pm2,\pm3\}$,
see for example \cite[Section 5]{han-exp}.

In this subsection, we write the elements of $G_{2}$ with respect
to the basis of $V_{7}$ since $G_{2}$ can be viewed as a Lie subalgebra
of $\mathfrak{gl}(V_{7})$. In this case $G_{2}=\langle x_{i},y_{i},h_{1},h_{2}:i=1,\dots,6\rangle$
where 
\[
x_{1}=\begin{pmatrix}0 & -1 & 0 & 0 & 0 & 0 & 0\\
0 & 0 & 0 & 0 & 0 & 0 & 0\\
0 & 0 & 0 & 1 & 0 & 0 & 0\\
0 & 0 & 0 & 0 & -2 & 0 & 0\\
0 & 0 & 0 & 0 & 0 & 0 & 0\\
0 & 0 & 0 & 0 & 0 & 0 & 1\\
0 & 0 & 0 & 0 & 0 & 0 & 0
\end{pmatrix},\ x_{2}=\begin{pmatrix}0 & 0 & 0 & 0 & 0 & 0 & 0\\
0 & 0 & 1 & 0 & 0 & 0 & 0\\
0 & 0 & 0 & 0 & 0 & 0 & 0\\
0 & 0 & 0 & 0 & 0 & 0 & 0\\
0 & 0 & 0 & 0 & 0 & -1 & 0\\
0 & 0 & 0 & 0 & 0 & 0 & 0\\
0 & 0 & 0 & 0 & 0 & 0 & 0
\end{pmatrix},
\]
\[
y_{1}=\begin{pmatrix}0 & 0 & 0 & 0 & 0 & 0 & 0\\
-1 & 0 & 0 & 0 & 0 & 0 & 0\\
0 & 0 & 0 & 0 & 0 & 0 & 0\\
0 & 0 & 2 & 0 & 0 & 0 & 0\\
0 & 0 & 0 & -1 & 0 & 0 & 0\\
0 & 0 & 0 & 0 & 0 & 0 & 0\\
0 & 0 & 0 & 0 & 0 & 1 & 0
\end{pmatrix},\ y_{2}=\begin{pmatrix}0 & 0 & 0 & 0 & 0 & 0 & 0\\
0 & 0 & 0 & 0 & 0 & 0 & 0\\
0 & 1 & 0 & 0 & 0 & 0 & 0\\
0 & 0 & 0 & 0 & 0 & 0 & 0\\
0 & 0 & 0 & 0 & 0 & 0 & 0\\
0 & 0 & 0 & 0 & -1 & 0 & 0\\
0 & 0 & 0 & 0 & 0 & 0 & 0
\end{pmatrix},
\]
 $h_{1}=\text{diag}(1,-1,2,0,-2,1,-1)$, $h_{2}=\text{diag}(0,1,-1,0,1,-1,0)$,
and $x_{3}=[x_{1},x_{2}]$, $x_{4}=[x_{1},x_{3}]$, $x_{5}=[x_{1},x_{4}]$,
$x_{6}=[x_{5},x_{2}]$, $y_{3}=[y_{1},y_{2}]$, $y_{4}=[y_{1},y_{3}]$,
$y_{5}=[y_{1},y_{4}]$, $y_{6}=[y_{5},y_{2}]$. 

We know that the map $\mathfrak{g}_{\bar{0}}\times\mathfrak{g}_{\bar{0}}\rightarrow\mathfrak{g}_{\bar{0}}$
is a Lie bracket and the Lie superbracket $[\cdotp,\cdotp]:\mathfrak{g}_{\bar{0}}\times\mathfrak{g}_{\bar{1}}\rightarrow\mathfrak{g}_{\bar{1}}$
is defined by $[x+y,v_{i}\otimes e_{j}]=xv_{i}\otimes e_{j}+v_{i}\otimes ye_{j}$
for $x\in\mathfrak{sl}_{2}(\mathbb{K})$, $y\in G_{2}$ and $v_{i}\otimes e_{j}\in\mathfrak{g}_{\bar{1}}$.
The Lie superbracket $[\cdotp,\cdotp]:\mathfrak{g}_{\bar{1}}\times\mathfrak{g}_{\bar{1}}\rightarrow\mathfrak{g}_{\bar{0}}$
over a field of zero characteristic is calculated explicitly in \cite[Subsection 5.1]{han-exp}.
Note that the structure remains the same when we work in the field
$\mathbb{K}$ with $\mathrm{char}(\mathbb{K})=p>3$ and we adopt it
to calculate a basis for $\mathfrak{g}^{e}$ and $\mathfrak{z}(\mathfrak{g}^{e})$.

It is clear that representatives of nilpotent orbits in $\mathfrak{sl}_{2}(\mathbb{K})$
up to conjugation by $\mathrm{SL}_{2}(\mathbb{K})$ are $0$ and $E$.
Based on \cite[Section 11]{Lawther2008}, representatives of nilpotent
orbits in Lie algebra $G_{2}$ up to conjugation by Lie group $G_{2}$
are $0$, $x_{1}$, $x_{2}$, $x_{2}+x_{5}$ and $x_{1}+x_{2}$. Representatives
of nilpotent orbits $e\in G(3)$ are listed in Tables \ref{tab:G(3)-1}--\ref{tab:G(3)-2}.
In these tables we also give the cocharacter $\tau$ associated to
$e$ and basis elements for $\mathfrak{g}^{e}$ and $\mathfrak{z}\left(\mathfrak{g}^{e}\right)$.
Note that in the second column of Tables \ref{tab:G(3)-1} and \ref{tab:G(3)-2},
we denote by $\alpha_{0}=2\delta$, $\alpha_{1}=\varepsilon_{1}$,
$\alpha_{2}=\varepsilon_{2}-\varepsilon_{1}$ where $2\delta,\varepsilon_{1},\varepsilon_{2}-\varepsilon_{1}$
are even roots in the root system for $\mathfrak{g}$ as given in
Section \ref{sec:highest root}. After Table \ref{tab:G(3)-2}, we
take the case when $e=x_{2}$ as an example to show our calculation
explicitly. The other cases are obtained using a similar approach.

\begin{table}[H]
\noindent \begin{centering}
\begin{tabular}{|>{\centering}p{2.4cm}|>{\centering}p{3.2cm}|>{\centering}p{4.7cm}|c|c|}
\hline 
Representatives of nilpotent orbits $e\in\mathfrak{g}_{\bar{0}}$ & Cocharacter $\tau$ & $\mathfrak{g}^{e}$ & $\dim\mathfrak{g}^{e}$ & $\mathfrak{z}\left(\mathfrak{g}^{e}\right)$\tabularnewline
\hline 
\hline 
$0$ & $0$ & $\mathfrak{g}$ & $31$ & $0$\tabularnewline
\hline 
$x_{2}$ & $h_{\alpha_{2}}(t)$ & $\langle E,H,F,2h_{1}+3h_{2},x_{2},x_{3},x_{6},y_{1},y_{5},x_{4},y_{4}\rangle\oplus\langle v_{i}\otimes e_{j}:i=\pm1,j=0,-1,2,\pm3\rangle$ & $21$ & $\langle e\rangle$\tabularnewline
\hline 
$x_{1}$ & $h_{\alpha_{1}}(t)$ & $\langle E,H,F,x_{5},y_{2},x_{1},x_{6},y_{6},h_{1}+2h_{2}\rangle\oplus\langle v_{i}\otimes e_{j}:i=\pm1,j=-2,1,3\rangle$ & $15$ & $\langle e\rangle$\tabularnewline
\hline 
$x_{2}+x_{5}$ & $h_{\alpha_{1}}(t^{2})h_{\alpha_{2}}(t^{4})$ & $\langle E,H,F,x_{6},x_{2}+x_{5},x_{3},x_{4}\rangle\oplus\langle v_{i}\otimes e_{j}:i=\pm1,j=0,2,3\rangle$ & $13$ & $\langle e,x_{6}\rangle$\tabularnewline
\hline 
$x_{1}+x_{2}$ & $h_{\alpha_{1}}(t^{6})h_{\alpha_{2}}(t^{10})$ & $\langle E,H,F,x_{6},x_{1}+x_{2}\rangle\oplus\langle v_{i}\otimes e_{3}:i=\pm1\rangle$ & $7$ & $\langle e,x_{6}\rangle$\tabularnewline
\hline 
\end{tabular}\caption{\label{tab:G(3)-1}Cocharacters and $\mathfrak{g}^{e}$ for $\mathfrak{g}=G(3)$}
\par\end{centering}
\end{table}

\begin{table}[H]
\noindent \begin{centering}
\begin{tabular}{|>{\centering}p{2.4cm}|>{\centering}p{3.2cm}|>{\centering}p{4.7cm}|c|>{\centering}p{1.5cm}|}
\hline 
Representatives of nilpotent orbits $e\in\mathfrak{g}_{\bar{0}}$ & Cocharacter $\tau$ & $\mathfrak{g}^{e}$ & $\dim\mathfrak{g}^{e}$ & $\mathfrak{z}\left(\mathfrak{g}^{e}\right)$\tabularnewline
\hline 
\hline 
$E$ & $h_{\alpha_{0}}(t)$ & $\langle E\rangle\oplus G_{2}\oplus\langle v_{1}\otimes e_{j}:j=0,\pm1,\pm2,\pm3\rangle$ & $22$ & $\langle e\rangle$\tabularnewline
\hline 
$E+x_{2}$ & $h_{\alpha_{0}}(t)h_{\alpha_{2}}(t)$ & $\langle E\rangle\oplus\langle2h_{1}+3h_{2},x_{2},x_{3},x_{6},y_{1},y_{5},x_{4},y_{4}\rangle\oplus\langle v_{1}\otimes e_{2},v_{1}\otimes e_{-1},v_{1}\otimes e_{\pm3},v_{1}\otimes e_{0},v_{1}\otimes e_{-2}+v_{-1}\otimes e_{-1},v_{1}\otimes e_{1}-v_{-1}\otimes e_{2}\rangle$ & $16$ & $\langle e\rangle$\tabularnewline
\hline 
$E+x_{1}$ & $h_{\alpha_{0}}(t)h_{\alpha_{1}}(t)$ & $\langle E\rangle\oplus\langle x_{5},y_{2},x_{1},x_{6},y_{6},h_{1}+2h_{2}\rangle\oplus\langle v_{1}\otimes e_{1},v_{1}\otimes e_{3},v_{1}\otimes e_{-2},v_{1}\otimes e_{0}-v_{-1}\otimes e_{1},v_{1}\otimes e_{-3}-v_{-1}\otimes e_{-2},v_{1}\otimes e_{2}+v_{-1}\otimes e_{3}\rangle$ & $13$ & $\langle e\rangle$\tabularnewline
\hline 
$E+(x_{2}+x_{5})$ & $h_{\alpha_{0}}(t)h_{\alpha_{1}}(t^{2})h_{\alpha_{2}}(t^{4})$ & $\langle E\rangle\oplus\langle x_{6},x_{2}+x_{5},x_{3},x_{4}\rangle\oplus\langle v_{1}\otimes e_{3},v_{1}\otimes e_{2},v_{1}\otimes e_{0},v_{1}\otimes e_{1}-v_{-1}\otimes e_{2},6v_{-1}\otimes e_{3}-v_{1}\otimes e_{-1}\rangle$ & $10$ & $\langle e,x_{6}\rangle$\tabularnewline
\hline 
$E+(x_{1}+x_{2})$ & $h_{\alpha_{0}}(t)h_{\alpha_{1}}(t^{6})h_{\alpha_{2}}(t^{10})$ & $\langle E\rangle\oplus\langle x_{6},x_{1}+x_{2}\rangle\oplus\langle v_{1}\otimes e_{3},v_{-1}\otimes e_{3}+v_{1}\otimes e_{2}\rangle$ & $5$ & $\langle e,x_{6},v_{1}\otimes e_{3}\rangle$\tabularnewline
\hline 
\end{tabular}
\par\end{centering}
\caption{\label{tab:G(3)-2}Cocharacters and $\mathfrak{g}^{e}$ for $\mathfrak{g}=G(3)$
(continued)}

\end{table}

When $e=x_{2}$, we first compute that the cocharacter associated
to $e$ is $\tau(t)=h_{\alpha_{2}}(t)$. According to \cite[Section 5]{Lawther2008}
and \cite[Section 5.3]{han-exp}, we have that $\mathfrak{g}_{\bar{0}}^{e}=\mathfrak{sl}_{2}(\mathbb{K})\oplus\langle2h_{1}+3h_{2},x_{2},x_{3},x_{6},y_{1},y_{5},x_{4},y_{4}\rangle$.

In order to determine $\mathfrak{g}_{\bar{1}}^{e}$, we work out the
$\tau$-grading on basis elements of $\mathfrak{g}_{\bar{1}}$ in
the table below.

\begin{table}[H]
\noindent \begin{centering}
\begin{tabular}{|c|c|}
\hline 
$\tau$-grading & basis elements of $\mathfrak{g}_{\bar{1}}$\tabularnewline
\hline 
\hline 
$1$ & $v_{i}\otimes e_{2},v_{i}\otimes e_{-1}$, $i=\pm1$\tabularnewline
\hline 
$0$ & $v_{i}\otimes e_{\pm3},v_{i}\otimes e_{0}$, $i=\pm1$\tabularnewline
\hline 
$-1$ & $v_{i}\otimes e_{-2},v_{i}\otimes e_{1}$, $i=\pm1$\tabularnewline
\hline 
\end{tabular}
\par\end{centering}
\caption{The $\tau$-grading on basis elements of $\mathfrak{g}_{\bar{1}}$
when $e=x_{2}$}
\end{table}
 From the above table we know that $\mathfrak{g}_{\bar{1}}^{e}=\mathfrak{g}_{\bar{1}}^{e}(-1)\oplus\mathfrak{g}_{\bar{1}}^{e}(0)\oplus\mathfrak{g}_{\bar{1}}^{e}(1)$.
Since $[e,\mathfrak{g}_{\bar{1}}^{e}(1)]\subseteq\mathfrak{g}_{\bar{1}}^{e}(3)$
and $\mathfrak{g}_{\bar{1}}^{e}(3)=0$, we get that $\mathfrak{g}_{\bar{1}}^{e}(1)=\langle v_{i}\otimes e_{2},v_{i}\otimes e_{-1}:i=\pm1\rangle$.
Similarly, we have that $\mathfrak{g}_{\bar{1}}^{e}(0)=\langle v_{i}\otimes e_{\pm3},v_{i}\otimes e_{0}:i=\pm1\rangle$.
Then let $x=\sum a_{i,-2}v_{i}\otimes e_{-2}+\sum a_{i,1}v_{i}\otimes e_{1}\in\mathfrak{g}_{\bar{1}}^{e}(-1)$,
we have that $[e,x]=-\sum a_{i,-2}v_{i}\otimes e_{-1}+a_{i,1}v_{i}\otimes e_{2}$.
This equals to zero implies that $a_{i,-2}=a_{i,1}=0$, thus $\mathfrak{g}_{\bar{1}}^{e}(-1)=0$.
Hence, we deduce that $\mathfrak{g}_{\bar{1}}^{e}=\langle v_{i}\otimes e_{j}:i=\pm1,j=0,-1,2,\pm3\rangle$.

Next we look for $\mathfrak{z}(\mathfrak{g}^{e})=\mathfrak{z}_{\bar{0}}\oplus\mathfrak{z}_{\bar{1}}$.
Note that $\mathfrak{z}_{\bar{0}}=\mathfrak{z}_{\bar{0}}(2)\oplus\mathfrak{z}_{\bar{0}}(1)\oplus\mathfrak{z}_{\bar{0}}(0)$.
We know that there is no centre in $\mathfrak{sl}_{2}(\mathbb{K})$.
Let $y=b_{1}y_{1}+b_{3}x_{3}+b_{5}y_{5}+b_{6}x_{6}\in\mathfrak{z}_{\bar{0}}(1)$,
we compute $[2h_{1}+3h_{2},y]=-b_{1}y_{1}+b_{3}x_{3}-3b_{5}y_{5}+3b_{6}x_{6}$
for some $b_{i}\in\mathbb{K}$. This equals to zero if and only if
$b_{1}=b_{3}=b_{5}=b_{6}=0$. Thus $\mathfrak{z}_{\bar{0}}(1)=0$.
Using a similar method we have that $\mathfrak{z}_{\bar{0}}(0)=0$
and $\mathfrak{z}_{\bar{0}}(2)=\langle e\rangle$. For an element
$z\in\mathfrak{g}_{\bar{1}}^{e}$ such that $z=\sum_{i}a_{i,2}v_{i}\otimes e_{2}+\sum_{i}a_{i,-1}v_{i}\otimes e_{-1}+\sum_{i}a_{1,3}v_{i}\otimes e_{3}+\sum_{i}a_{1,-3}v_{i}\otimes e_{-3}+\sum_{i}a_{i,0}v_{i}\otimes e_{0}$,
by calculating $[H,z]$ we obtain that $\mathfrak{z}_{\bar{1}}=0$.
Therefore, we deduce that $\mathfrak{z}(\mathfrak{g}^{e})=\langle e\rangle$.

\subsection{Centralizer of even nilpotent elements in $F(4)$ and its centre\label{subsec:F(4)}}

Let $\mathfrak{g}=\mathfrak{g}_{\bar{0}}\oplus\mathfrak{g}_{\bar{1}}=F(4)$.
By definition, $\mathfrak{g}_{\bar{0}}=\mathfrak{sl}_{2}(\mathbb{K})\oplus\mathfrak{so}_{7}(\mathbb{K})$
and $\mathfrak{g}_{\bar{1}}=V_{2}\otimes V_{8}$ where $V_{2}=\langle v_{1}=(1,0)^{t},v_{-1}=(0,1)^{t}\rangle$
is the standard two-dimensional representation for $\mathfrak{sl}_{2}(\mathbb{K})$
and $V_{8}$ is the spin representation for $\mathfrak{so}_{7}(\mathbb{K})$.
Let $V_{\mathfrak{so}}$ be the standard $7$-dimensional module of
$\mathfrak{so}_{7}(\mathbb{K})$ such that $\mathfrak{so}_{7}(\mathbb{K})=\mathfrak{so}(V_{\mathfrak{so}})$
and let $\beta$ be the symmetric bilinear form on $V_{\mathfrak{so}}$.
We adopt the notation from \cite[Section 6.1]{han-exp} for a basis
for $V_{\mathfrak{so}}$. Note that there exists a decomposition of
$V_{\mathfrak{so}}$ such that $V_{\mathfrak{so}}=W\oplus\langle e_{0}\rangle\oplus W^{*}$
where $W$, $W^{*}$ is a pair of dual maximal isotropic subspaces
of $V_{\mathfrak{so}}$ corresponding to $\beta$ and $W=\langle e_{1},e_{2},e_{3}\rangle$,
$W^{*}=\langle e_{-1},e_{-2},e_{-3}\rangle$. The symmetric form $\beta$
on the chosen basis elements is given by $\beta(e_{0},e_{0})=2$,
$\beta(e_{0},W)=\beta(e_{0},W^{*})=0$ and $\beta(e_{i},e_{-j})=\delta_{ij}$
for $i,j\neq0$. Denote by $e_{i,j}$ the elementary transformation
which sends $e_{i}$ to $e_{j}$ and the other basis vectors to $0$.
Then the following elements 
\begin{equation}
R_{i,-j}=e_{i,j}-e_{-j,-i}\text{ and }R_{i,0}=2e_{i,0}-e_{0,-i}\text{ for }i,j\in\{\pm1,\pm2,\pm3\}\label{eq:Ri,j}
\end{equation}
 form a basis for $\mathfrak{so}_{7}(\mathbb{K})$, see \cite[Lemma 6.2.1]{Goodman}.
Next consider a $1$-dimensional vector space $\langle s\rangle$
such that $e_{-i}s=0$ for $i=1,2,3$ and $e_{0}s=s$. Then a basis
for $V_{8}$ can be written as $\{s$, $e_{1}s$, $e_{2}s$, $e_{3}s$,
$e_{1}e_{2}s$, $e_{1}e_{3}s$, $e_{2}e_{3}s$, $e_{1}e_{2}e_{3}s\}$,
see \cite[Subsection 6.2]{han-exp}.

According to \cite[Section 1.6]{Janzten}, nilpotent orbits in $\mathfrak{so}_{7}(\mathbb{K})$
are parameterized by partition $\lambda$ such that $\lambda\in\{(7),(5,1^{2}),(3^{2},1),(3,1^{4}),(2^{2},3),(1^{7})\}$.
In Table \ref{tab:F(4)-centralizer}, we list representatives of nilpotent
orbits $e=e_{1}+e_{2}\in\mathfrak{sl}_{2}(\mathbb{K})\oplus\mathfrak{so}_{7}(\mathbb{K})$
up to conjugation by $\mathrm{SL}_{2}(\mathbb{K})\times\mathrm{SO}_{7}(\mathbb{K})$
such that $e_{2}$ are expressed in terms of $R_{i,j}$ in (\ref{eq:Ri,j}).
Basis elements of $\mathfrak{so}_{7}(\mathbb{C})^{e_{2}}$ are given
in \cite[Table 10]{han-exp}, note that the basis of $\mathfrak{so}_{7}(\mathbb{C})^{e_{2}}$
viewed in $\mathfrak{so}_{7}(\mathbb{K})^{e_{2}}$ is also a basis
for $\mathfrak{so}_{7}(\mathbb{K})^{e_{2}}$ according to \cite[Section 3.2]{Janzten}.
We also give the cocharacter $\tau$ associated to $e$ and recall
basis elements for $\mathfrak{g}_{\bar{0}}^{e}$ in the table below.
Note that in the second column of Table \ref{tab:F(4)-centralizer},
we denote by $\alpha_{0}=\delta$, $\alpha_{1}=\varepsilon_{1}-\varepsilon_{2}$,
$\alpha_{2}=\varepsilon_{2}-\varepsilon_{3}$ and $\alpha_{3}=\varepsilon_{3}$
where $\delta$, $\varepsilon_{1}-\varepsilon_{2}$, $\varepsilon_{2}-\varepsilon_{3}$,
$\varepsilon_{3}$ are even roots in the root system for $\mathfrak{g}$
as given in Section \ref{sec:highest root}.
\begin{singlespace}
\noindent \begin{center}
\begin{longtable}[c]{|>{\centering}m{2.7cm}|>{\centering}m{4.2cm}|>{\centering}m{6cm}|>{\centering}m{0.9cm}|}
\caption{\label{tab:F(4)-centralizer}Cocharacters and $\mathfrak{g}_{\bar{0}}^{e}$
for $\mathfrak{g}=F(4)$}
\tabularnewline
\endfirsthead
\hline 
Representatives of nilpotent orbits $e\in\mathfrak{g}_{\bar{0}}$ & Cocharacter $\tau$ & $\mathfrak{g}_{\bar{0}}^{e}$ & $\dim\mathfrak{g}_{\bar{0}}^{e}$\tabularnewline
\hline 
\hline 
$e_{(7)}=R_{1,-2}+R_{2,-3}+R_{3,0}$ & $h_{\alpha_{1}}(t^{6})h_{\alpha_{2}}(t^{10})h_{\alpha_{3}}(t^{6})$ & $\mathfrak{sl}_{2}(\mathbb{K})\oplus\langle e,R_{1,0}-2R_{2,3},R_{1,2}\rangle$ & $6$\tabularnewline
\hline 
$e_{(5,1^{2})}=R_{1,-2}+R_{2,0}$ & $h_{\alpha_{1}}(t^{4})h_{\alpha_{2}}(t^{6})h_{\alpha_{3}}(t^{3})$ & $\mathfrak{sl}_{2}(\mathbb{K})\oplus\langle R_{3,-3},e,R_{1,-3},R_{1,3},R_{1,2}\rangle$ & $8$\tabularnewline
\hline 
$e_{(3^{2},1)}=R_{1,-3}+R_{2,3}$ & $h_{\alpha_{1}}(t^{2})h_{\alpha_{2}}(t^{4})h_{\alpha_{3}}(t^{2})$ & $\mathfrak{sl}_{2}(\mathbb{K})\oplus\langle R_{1,-1}-R_{2,-2}+R_{3,-3},$

$e,R_{2,-3},R_{2,0},R_{1,0},R_{1,3},R_{1,2}\rangle$ & $10$\tabularnewline
\hline 
$e_{(3,2^{2})}=R_{1,0}+R_{2,3}$ & $h_{\alpha_{1}}(t^{2})h_{\alpha_{2}}(t^{3})h_{\alpha_{3}}(t^{2})$ & $\mathfrak{sl}_{2}(\mathbb{K})\oplus\langle R_{2,-2}-R_{3,-3},R_{2,-3},$

$R_{3,-2},2R_{1,-3}+R_{2,0},-2R_{1,-2}+R_{3,0},$

$R_{1,0},R_{2,3},R_{1,3},R_{1,2}\rangle$ & $12$\tabularnewline
\hline 
$e_{(3,1^{4})}=R_{1,0}$ & $h_{\alpha_{1}}(t^{2})h_{\alpha_{2}}(t^{2})h_{\alpha_{3}}(t)$ & $\mathfrak{sl}_{2}(\mathbb{K})\oplus\langle R_{2,-2},R_{3,-3},R_{2,3},R_{2,-3},$

$R_{-3,-2},R_{3,-2},e,R_{1,2},R_{1,3},R_{1,-3},R_{1,-2}\rangle$ & $14$\tabularnewline
\hline 
$e_{(2^{2},1^{3})}=R_{1,2}$ & $h_{\alpha_{1}}(t)h_{\alpha_{2}}(t^{2})h_{\alpha_{3}}(t)$ & $\mathfrak{sl}_{2}(\mathbb{K})\oplus\langle R_{1,-2},R_{1,-1}-R_{2,-2},R_{3,-3},$

$R_{2,-1},R_{-3,0},R_{3,0},R_{1,3},R_{1,0},R_{2,-3},$

$R_{2,3},R_{1,-3},R_{2,0},e\rangle$ & $16$\tabularnewline
\hline 
$e_{(1^{7})}=0$ & $0$ & $\mathfrak{sl}_{2}(\mathbb{K})\oplus\mathfrak{so}_{7}(\mathbb{K})$ & $24$\tabularnewline
\hline 
$E+e_{(7)}$ & $h_{\alpha_{0}}(t)h_{\alpha_{1}}(t^{6})h_{\alpha_{2}}(t^{10})h_{\alpha_{3}}(t^{6})$ & $\langle E,e_{(7)},R_{1,0}-2R_{2,3},R_{1,2}\rangle$ & $4$\tabularnewline
\hline 
$E+e_{(5,1^{2})}$ & $h_{\alpha_{0}}(t)h_{\alpha_{1}}(t^{4})h_{\alpha_{2}}(t^{6})h_{\alpha_{3}}(t^{3})$ & $\langle E,R_{3,-3},e_{(5,1^{2})},R_{1,-3},R_{1,3},R_{1,2}\rangle$ & $6$\tabularnewline
\hline 
$E+e_{(3^{2},1)}$ & $h_{\alpha_{0}}(t)h_{\alpha_{1}}(t^{2})h_{\alpha_{2}}(t^{4})h_{\alpha_{3}}(t^{2})$ & $\langle E,R_{1,-1}-R_{2,-2}+R_{3,-3},e_{(3^{2},1)},$

$R_{2,-3},R_{2,0},R_{1,0},R_{1,3},R_{1,2}\rangle$ & $8$\tabularnewline
\hline 
$E+e_{(3,2^{2})}$ & $h_{\alpha_{0}}(t)h_{\alpha_{1}}(t^{2})h_{\alpha_{2}}(t^{3})h_{\alpha_{3}}(t^{2})$ & $\langle E,R_{2,-2}-R_{3,-3},R_{2,-3},R_{3,-2},$

$2R_{1,-3}+R_{2,0},-2R_{1,-2}+R_{3,0},$

$R_{1,0},R_{2,3},R_{1,3},R_{1,2}\rangle$ & $10$\tabularnewline
\hline 
$E+e_{(3,1^{4})}$ & $h_{\alpha_{0}}(t)h_{\alpha_{1}}(t^{2})h_{\alpha_{2}}(t^{2})h_{\alpha_{3}}(t)$ & $\langle E,R_{2,-2},R_{3,-3},R_{2,3},R_{2,-3},R_{-3,-2},$

$R_{3,-2},e_{(3,1^{4})},R_{1,2},R_{1,3},R_{1,-3},R_{1,-2}\rangle$ & $12$\tabularnewline
\hline 
$E+e_{(2^{2},1^{3})}$ & $h_{\alpha_{0}}(t)h_{\alpha_{1}}(t)h_{\alpha_{2}}(t^{2})h_{\alpha_{3}}(t)$ & $\langle E,R_{1,-2},R_{1,-1}-R_{2,-2},R_{3,-3},R_{2,-1},$

$R_{-3,0},R_{3,0},R_{1,3},R_{1,0},R_{2,-3},R_{2,3},$

$R_{1,-3},R_{2,0},e_{(2^{2},1^{3})}\rangle$ & $14$\tabularnewline
\hline 
$E$ & $h_{\alpha_{0}}(t)$ & $\langle E\rangle\oplus\mathfrak{so}_{7}(\mathbb{K})$ & $22$\tabularnewline
\hline 
\end{longtable}
\par\end{center}
\end{singlespace}

We then give basis elements for $\mathfrak{g}_{\bar{1}}^{e}$ and
$\mathfrak{z}(\mathfrak{g}^{e})$ in Table \ref{tab:z(g^e)-F(4)}.
\begin{singlespace}
\noindent \begin{center}
\begin{longtable}[c]{|>{\centering}m{2.3cm}|>{\centering}m{6.8cm}|>{\centering}m{1cm}|>{\centering}m{2cm}|}
\caption{\label{tab:z(g^e)-F(4)}$\mathfrak{g}_{\bar{1}}^{e}$ and $\mathfrak{z}(\mathfrak{g}^{e})$
for $\mathfrak{g}=F(4)$}
\tabularnewline
\endfirsthead
\hline 
Representatives of nilpotent orbits $e\in\mathfrak{g}_{\bar{0}}$ & $\mathfrak{g}_{\bar{1}}^{e}$ & $\dim\mathfrak{g}_{\bar{1}}^{e}$ & $\mathfrak{z}(\mathfrak{g}^{e})$\tabularnewline
\hline 
\hline 
$e_{(7)}$ & $\langle v_{i}\otimes e_{1}e_{2}e_{3}s,v_{i}\otimes e_{1}s-v_{i}\otimes e_{2}e_{3}s:i=\pm1\rangle$ & $4$ & $\langle e,R_{1,2}\rangle$\tabularnewline
\hline 
$e_{(5,1^{2})}$ & $\langle v_{i}\otimes e_{1}e_{2}e_{3}s,v_{i}\otimes e_{1}e_{2}s:i=\pm1\rangle$ & $4$ & $\langle e,R_{1,2}\rangle$\tabularnewline
\hline 
$e_{(3^{2},1)}$ & $\langle v_{i}\otimes e_{1}e_{2}s,v_{i}\otimes e_{1}e_{2}e_{3}s,v_{i}\otimes e_{1}e_{3}s,$

$v_{i}\otimes e_{2}s:i=\pm1\rangle$ & $8$ & $\langle e,R_{1,2}\rangle$\tabularnewline
\hline 
$e_{(3,2^{2})}$ & $\langle v_{i}\otimes e_{1}e_{2}e_{3}s,v_{i}\otimes e_{1}e_{2}s,v_{i}\otimes e_{1}e_{3}s,$

$v_{i}\otimes e_{1}s-v_{i}\otimes e_{2}e_{3}s:i=\pm1\rangle$ & $8$ & $\langle e\rangle$\tabularnewline
\hline 
$e_{(3,1^{4})}$ & $\langle v_{i}\otimes e_{1}e_{2}e_{3}s,v_{i}\otimes e_{1}e_{2}s,v_{i}\otimes e_{1}e_{3}s,$

$v_{i}\otimes e_{1}s:i=\pm1\rangle$ & $8$ & $\langle e\rangle$\tabularnewline
\hline 
$e_{(2^{2},1^{3})}$ & $\langle v_{i}\otimes e_{1}e_{2}e_{3}s,v_{i}\otimes e_{1}e_{2}s,v_{i}\otimes e_{1}e_{3}s,$

$v_{i}\otimes e_{1}s,v_{i}\otimes e_{2}s,v_{i}\otimes e_{2}e_{3}s:i=\pm1\rangle$ & $12$ & $\langle e\rangle$\tabularnewline
\hline 
$e_{(1^{7})}$ & $\mathfrak{g}_{\bar{1}}$ & $16$ & $\langle e\rangle=0$\tabularnewline
\hline 
$E+e_{(7)}$ & $\langle v_{1}\otimes e_{1}e_{2}e_{3}s,v_{1}\otimes e_{1}e_{2}s-v_{-1}\otimes e_{1}e_{2}e_{3}s,$

$v_{1}\otimes e_{1}s-v_{1}\otimes e_{2}e_{3}s\rangle$ & $3$ & $\langle e,R_{1,2},$

$v_{1}\otimes e_{1}e_{2}e_{3}s\rangle$\tabularnewline
\hline 
$E+e_{(5,1^{2})}$ & $\langle v_{1}\otimes e_{1}e_{2}e_{3}s,v_{1}\otimes e_{1}e_{2}s,v_{1}\otimes e_{1}s$

$-v_{-1}\otimes e_{1}e_{2}s,v_{1}\otimes e_{1}e_{3}s+v_{-1}\otimes e_{1}e_{2}e_{3}s\rangle$ & $4$ & $\langle e,R_{1,2}\rangle$\tabularnewline
\hline 
$E+e_{(3^{2},1)}$ & $\langle v_{1}\otimes e_{1}e_{2}s,v_{1}\otimes e_{1}e_{2}e_{3}s,v_{1}\otimes e_{2}s,$

$v_{1}\otimes e_{1}s-v_{-1}\otimes e_{1}e_{2}e_{3}s,v_{1}\otimes e_{1}e_{3}s,$

$v_{-1}\otimes e_{1}e_{2}s+v_{1}\otimes e_{2}e_{3}s\rangle$ & $6$ & $\langle e,R_{1,2}\rangle$\tabularnewline
\hline 
$E+e_{(3,2^{2})}$ & $\langle v_{1}\otimes e_{1}e_{2}e_{3}s,v_{1}\otimes e_{1}e_{2}s,v_{1}\otimes e_{1}e_{3}s,$

$v_{1}\otimes e_{1}s-v_{-1}\otimes e_{1}e_{2}e_{3}s,v_{1}\otimes e_{2}e_{3}s$

$-v_{-1}\otimes e_{1}e_{2}e_{3}s,v_{1}\otimes e_{3}s+v_{-1}\otimes e_{1}e_{3}s,$

$v_{1}\otimes e_{2}s+v_{-1}\otimes e_{1}e_{2}s\rangle$ & $7$ & $\langle e\rangle$\tabularnewline
\hline 
$E+e_{(3,1^{4})}$ & $\langle v_{1}\otimes e_{1}s,v_{1}\otimes e_{1}e_{3}s,v_{1}\otimes e_{1}e_{2}e_{3}s,$

$v_{1}\otimes e_{1}e_{2}s,v_{1}\otimes s-v_{-1}\otimes e_{1}s,v_{1}\otimes e_{2}s$

$+v_{-1}\otimes e_{1}e_{2}s,v_{1}\otimes e_{3}s+v_{-1}\otimes e_{1}e_{3}s,$ 

$v_{1}\otimes e_{2}e_{3}s-v_{-1}\otimes e_{1}e_{2}e_{3}s,\rangle$ & $8$ & $\langle e\rangle$\tabularnewline
\hline 
$E+e_{(2^{2},1^{3})}$ & $\langle v_{1}\otimes e_{1}e_{2}e_{3}s,v_{1}\otimes e_{1}e_{2}s,v_{1}\otimes e_{1}s,$

$v_{1}\otimes e_{2}s,v_{1}\otimes e_{1}e_{3}s,v_{1}\otimes e_{2}e_{3}s,v_{1}\otimes s$

$-v_{-1}\otimes e_{1}e_{2}s,v_{1}\otimes e_{3}s-v_{-1}\otimes e_{1}e_{2}e_{3}s\rangle$ & $8$ & $\langle e\rangle$\tabularnewline
\hline 
$E$ & $\langle v_{1}\otimes s,v_{1}\otimes e_{1}s,v_{1}\otimes e_{2}s,v_{1}\otimes e_{3}s,v_{1}\otimes e_{1}e_{2}s,$

$v_{1}\otimes e_{1}e_{3}s,v_{1}\otimes e_{2}e_{3}s,v_{1}\otimes e_{1}e_{2}e_{3}s\rangle$ & $8$ & $\langle e\rangle$\tabularnewline
\hline 
\end{longtable}
\par\end{center}
\end{singlespace}

In the remaining part of this subsection, we give explicit calculations
on finding $\mathfrak{g}^{e}$ and $\mathfrak{z}(\mathfrak{g}^{e})$
for the case $e=e_{(5,1^{2})}$. All the other cases are obtained
using a similar approach.

Since $e=e_{(5,1^{2})}=R_{1,-2}+R_{2,0}\in\mathfrak{so}_{7}(\mathbb{K})$,
by calculating $[e,x]=0$ for any $x\in\mathfrak{so}_{7}(\mathbb{K})$
we obtain that $\mathfrak{so}_{7}(\mathbb{K})^{e}=\langle R_{3,-3},e,R_{1,-3},R_{1,3},R_{1,2}\rangle$
and thus $\mathfrak{g}_{\bar{0}}^{e}=\langle E,H,F\rangle\oplus\mathfrak{so}_{7}(\mathbb{K})^{e}$.

To compute $\mathfrak{g}_{\bar{1}}^{e}$, we first work out the $\tau$-grading
on basis elements of $\mathfrak{g}_{\bar{1}}$ in the table below.

\begin{table}[H]
\noindent \begin{centering}
\begin{tabular}{|c|c|}
\hline 
$\tau$-grading & basis elements of $\mathfrak{g}_{\bar{1}}$\tabularnewline
\hline 
\hline 
$3$ & $v_{i}\otimes e_{1}e_{2}e_{3}s,v_{i}\otimes e_{1}e_{2}s\text{ for }i=\pm1$\tabularnewline
\hline 
$1$ & $v_{i}\otimes e_{1}s,v_{i}\otimes e_{1}e_{3}s\text{ for }i=\pm1$\tabularnewline
\hline 
$-1$ & $v_{i}\otimes e_{2}s,v_{i}\otimes e_{2}e_{3}s\text{ for }i=\pm1$\tabularnewline
\hline 
$-3$ & $v_{i}\otimes s,v_{i}\otimes e_{3}s\text{ for }i=\pm1$\tabularnewline
\hline 
\end{tabular}
\par\end{centering}
\caption{The $\tau$-grading on basis elements of $\mathfrak{g}_{\bar{1}}$
when $e=e_{(5,1^{2})}$}

\end{table}

According to the table above, we have that $\mathfrak{g}_{\bar{1}}^{e}=\mathfrak{g}_{\bar{1}}^{e}(3)\oplus\mathfrak{g}_{\bar{1}}^{e}(1)\oplus\mathfrak{g}_{\bar{1}}^{e}(-1)\oplus\mathfrak{g}_{\bar{1}}^{e}(-3)$.
We get $\mathfrak{g}_{\bar{1}}^{e}(3)=\langle v_{i}\otimes e_{1}e_{2}e_{3}s,v_{i}\otimes e_{1}e_{2}s:i=\pm1\rangle$
directly from the grading. For any element $x=\sum_{i}a_{i,1}v_{i}\otimes e_{1}s+\sum_{i}a_{i,13}v_{i}\otimes e_{1}e_{3}s\in\mathfrak{g}_{\bar{1}}^{e}(1)$,
we have that $[e,x]=\sum_{i}a_{i,1}v_{i}\otimes e_{1}e_{2}s-\sum_{i}a_{i,13}v_{i}\otimes e_{1}e_{2}e_{3}s$.
Thus $[e,x]=0$ if and only if $a_{i,1}=a_{i,13}=0$. Hence, we know
that $\mathfrak{g}_{\bar{1}}^{e}(1)=0$. Applying a similar calculation
we obtain that $\mathfrak{g}_{\bar{1}}^{e}(-1)=\mathfrak{g}_{\bar{1}}^{e}(-3)=0$
and therefore $\mathfrak{g}_{\bar{1}}^{e}=\mathfrak{g}_{\bar{1}}^{e}(3)=\langle v_{i}\otimes e_{1}e_{2}e_{3}s,v_{i}\otimes e_{1}e_{2}s:i=\pm1\rangle$. 

Next we determine $\mathfrak{z}(\mathfrak{g}^{e})=\mathfrak{z}_{\bar{0}}\oplus\mathfrak{z}_{\bar{1}}$.
Clearly $e\in\mathfrak{z}_{\bar{0}}$. Let $y=b_{1}R_{3,-3}+b_{2}e+b_{3}R_{1,-3}+b_{4}R_{1,3}+b_{5}R_{1,2}$
be an element in $\mathfrak{z}_{\bar{0}}$. Then $[R_{3,-3},y]=-b_{3}R_{1,-3}+b_{4}R_{1,3}$,
this equals to zero implies that $b_{3}=b_{4}=0$. Similarly by computing
$[R_{1,3},y]$ we deduce that $b_{1}=0$. Note that $R_{1,2}\in\mathfrak{g}^{e}(6)$
and $\mathfrak{g}^{e}(6+j)=0$ for all $j$ such that $\mathfrak{g}^{e}(j)\neq0$.
Thus we deduce that $R_{1,2}\in\mathfrak{z}_{\bar{0}}$ and $\mathfrak{z}_{\bar{0}}=\langle e,R_{e_{1},e_{2}}\rangle$.
To determine $\mathfrak{z}_{\bar{1}}$, let $z=\sum_{i}a_{i,123}v_{i}\otimes e_{1}e_{2}e_{3}+\sum_{i}a_{i,12}v_{i}\otimes e_{1}e_{2}s\in\mathfrak{z}_{\bar{1}}$.
Then $[H,z]=\sum_{i}ia_{i,123}v_{i}\otimes e_{1}e_{2}e_{3}+\sum_{i}ia_{i,12}v_{i}\otimes e_{1}e_{2}s$,
this equals to zero if and only if $a_{i,123}=a_{i,12}=0$. Hence
we deduce that $\mathfrak{z}_{\bar{1}}=0$. Therefore, we conclude
that $\mathfrak{z}(\mathfrak{g}^{e})=\langle e,R_{e_{1},e_{2}}\rangle$.

\section{Reachability and the Panyushev property for exceptional Lie superalgebras\label{sec:Reachability}}

\noindent Let $\mathfrak{g}=\mathfrak{g}_{\bar{0}}\oplus\mathfrak{g}_{\bar{1}}$
be one of exceptional Lie superalgebras $D(2,1;\alpha)$, $G(3)$
and $F(4)$ over $\mathbb{K}$. With the structure of $\mathfrak{g}^{e}$
obtained in Section \ref{sec:exceptional}, we further investigate
some properties of nilpotent elements $e\in\mathfrak{g}_{\bar{0}}$.
Recall as in Introduction that $e$ is called \textit{reachable} if
$e\in[\mathfrak{g}^{e},\mathfrak{g}^{e}]$ and is called \textit{strongly
reachable} if $\mathfrak{g}^{e}=[\mathfrak{g}^{e},\mathfrak{g}^{e}]$.
The element $e$ is said to satisfy the \textit{Panyushev property}
if in the $\tau$-grading $\mathfrak{g}^{e}=\bigoplus_{j\geq0}\mathfrak{g}(j)$,
the subalgebra $\mathfrak{g}^{e}(\geq1)=\bigoplus_{j\geq1}\mathfrak{g}(j)$
is generated by $\mathfrak{g}^{e}(1)$. 

In this section, we adopt notation given in Section \ref{sec:exceptional}
for representatives of nilpotent orbits $e\in\mathfrak{g}_{\bar{0}}$
and basis elements of $\mathfrak{g}^{e}$. For each $e\in\mathfrak{g}_{\bar{0}}$,
we identify whether it is reachable, strongly reachable or satisfies
the Panyushev property in the following table. Note that the classification
is the same as the case of zero characteristic except for $\mathfrak{g}=G(3)$,
$e=E+x_{1}$, in which case we have that $e$ is reachable whenever
$p>5$.

\begin{table}[H]

\noindent \begin{centering}
\begin{tabular}{|>{\centering}p{3cm}|>{\centering}p{3cm}|>{\centering}p{2cm}|>{\centering}p{2cm}|>{\centering}p{3cm}|}
\hline 
Lie superalgebra $\mathfrak{g}$ & Representatives of nilpotent orbits $e\in\mathfrak{g}_{\bar{0}}$ & Reachable & Strongly reachable & Satisfying the Panyushev property\tabularnewline
\hline 
\hline 
\multirow{4}{3cm}{$D(2,1;\alpha)$} & $0$ & Yes & Yes & Yes\tabularnewline
\cline{2-5} \cline{3-5} \cline{4-5} \cline{5-5} 
 & $E_{1}$, $E_{2}$, $E_{3}$ & Yes & Yes & Yes\tabularnewline
\cline{2-5} \cline{3-5} \cline{4-5} \cline{5-5} 
 & $E_{1}+E_{2}$, $E_{2}+E_{3}$, $E_{1}+E_{3}$ & No & No & No\tabularnewline
\cline{2-5} \cline{3-5} \cline{4-5} \cline{5-5} 
 & $E_{1}+E_{2}+E_{3}$ & Yes & No & Yes\tabularnewline
\hline 
\end{tabular}\caption{\label{tab:Reachable D}Reachable, strongly reachable and Panyushev
elements in $D(2,1;\alpha)$}
\par\end{centering}
\end{table}

\begin{table}[H]
\noindent \begin{centering}
\begin{tabular}{|>{\centering}m{3cm}|>{\centering}m{3cm}|>{\centering}m{2cm}|>{\centering}m{2cm}|>{\centering}m{3cm}|}
\hline 
Lie superalgebra $\mathfrak{g}$ & Representatives of nilpotent orbits $e\in\mathfrak{g}_{\bar{0}}$ & Reachable & Strongly reachable & Satisfying the Panyushev property\tabularnewline
\hline 
\hline 
\multirow{10}{3cm}{$G(3)$} & $0$ & Yes & Yes & Yes\tabularnewline
\cline{2-5} \cline{3-5} \cline{4-5} \cline{5-5} 
 & $x_{2}$ & Yes & Yes & Yes\tabularnewline
\cline{2-5} \cline{3-5} \cline{4-5} \cline{5-5} 
 & $x_{1}$ & Yes & Yes & No\tabularnewline
\cline{2-5} \cline{3-5} \cline{4-5} \cline{5-5} 
 & $x_{2}+x_{5}$ & No & No & No\tabularnewline
\cline{2-5} \cline{3-5} \cline{4-5} \cline{5-5} 
 & $x_{1}+x_{2}$ & No & No & No\tabularnewline
\cline{2-5} \cline{3-5} \cline{4-5} \cline{5-5} 
 & $E$ & Yes & Yes & Yes\tabularnewline
\cline{2-5} \cline{3-5} \cline{4-5} \cline{5-5} 
 & $E+x_{2}$ & Yes & Yes & Yes\tabularnewline
\cline{2-5} \cline{3-5} \cline{4-5} \cline{5-5} 
 & $E+x_{1}$ & Yes except $p=5$ & No & No\tabularnewline
\cline{2-5} \cline{3-5} \cline{4-5} \cline{5-5} 
 & $E+(x_{2}+x_{5})$ & Yes & No & Yes\tabularnewline
\cline{2-5} \cline{3-5} \cline{4-5} \cline{5-5} 
 & $E+(x_{1}+x_{2})$ & No & No & No\tabularnewline
\hline 
\end{tabular}
\par\end{centering}
\caption{\label{tab:Reachable G}Reachable, strongly reachable and Panyushev
elements in $G(3)$}

\end{table}

\begin{table}[H]
\noindent \begin{centering}
\begin{tabular}{|>{\centering}m{3cm}|>{\centering}m{3cm}|>{\centering}m{2cm}|>{\centering}m{2cm}|>{\centering}m{3cm}|}
\hline 
Lie superalgebra $\mathfrak{g}$ & Representatives of nilpotent orbits $e\in\mathfrak{g}_{\bar{0}}$ & Reachable & Strongly reachable & Satisfying the Panyushev property\tabularnewline
\hline 
\hline 
\multirow{14}{3cm}{$F(4)$} & $e_{(7)}$ & No & No & No\tabularnewline
\cline{2-5} \cline{3-5} \cline{4-5} \cline{5-5} 
 & $e_{(5,1^{2})}$ & No & No & No\tabularnewline
\cline{2-5} \cline{3-5} \cline{4-5} \cline{5-5} 
 & $e_{(3^{2},1)}$ & No & No & No\tabularnewline
\cline{2-5} \cline{3-5} \cline{4-5} \cline{5-5} 
 & $e_{(3,2^{2})}$ & Yes & Yes & Yes\tabularnewline
\cline{2-5} \cline{3-5} \cline{4-5} \cline{5-5} 
 & $e_{(3,1^{4})}$ & Yes & Yes & Yes\tabularnewline
\cline{2-5} \cline{3-5} \cline{4-5} \cline{5-5} 
 & $e_{(2^{2},1^{3})}$ & Yes & Yes & Yes\tabularnewline
\cline{2-5} \cline{3-5} \cline{4-5} \cline{5-5} 
 & $e_{(1^{7})}$ & Yes & Yes & Yes\tabularnewline
\cline{2-5} \cline{3-5} \cline{4-5} \cline{5-5} 
 & $E+e_{(7)}$ & No & No & No\tabularnewline
\cline{2-5} \cline{3-5} \cline{4-5} \cline{5-5} 
 & $E+e_{(5,1^{2})}$ & No & No & No\tabularnewline
\cline{2-5} \cline{3-5} \cline{4-5} \cline{5-5} 
 & $E+e_{(3^{2},1)}$ & Yes & No & Yes\tabularnewline
\cline{2-5} \cline{3-5} \cline{4-5} \cline{5-5} 
 & $E+e_{(3.2^{2})}$ & Yes & No & Yes\tabularnewline
\cline{2-5} \cline{3-5} \cline{4-5} \cline{5-5} 
 & $E+e_{(3.1^{4})}$ & No & No & No\tabularnewline
\cline{2-5} \cline{3-5} \cline{4-5} \cline{5-5} 
 & $E+e_{(2^{2}.1^{3})}$ & Yes & Yes & Yes\tabularnewline
\cline{2-5} \cline{3-5} \cline{4-5} \cline{5-5} 
 & $E$ & Yes & Yes & Yes\tabularnewline
\hline 
\end{tabular}
\par\end{centering}
\caption{\label{tab:Reachable F}Reachable, strongly reachable and Panyushev
elements in $F(4)$}

\end{table}

In the remaining part of this section, we explain our calculations
explicitly for the following three cases: (1) $\mathfrak{g}=D(2,1;\alpha)$,
$e=E_{1}+E_{2}+E_{3}$; (2) $\mathfrak{g}=G(3)$, $e=E+x_{1}$; and
(3) $\mathfrak{g}=F(4)$, $e=E+e_{(3^{2},1)}$. Note that all other
cases can be done using a similar approach. 
\begin{itemize}
\item $\mathfrak{g}=D(2,1;\alpha)$, $e=E_{1}+E_{2}+E_{3}$
\end{itemize}
Based on Table \ref{tab:D(2,1)}, we have that $\mathfrak{g}^{e}=\mathfrak{g}^{e}(1)\oplus\mathfrak{g}^{e}(2)\oplus\mathfrak{g}^{e}(3)$
where $\mathfrak{g}^{e}(1)=\langle v_{1,1,-1}-v_{-1,1,1},v_{1,-1,1}-v_{-1,1,1}\rangle$,
$\mathfrak{g}^{e}(2)$=$\langle E_{1},E_{2},E_{3}\rangle$ and $\mathfrak{g}^{e}(3)=\langle v_{1,1,1}\rangle$. 

We first check whether $\mathfrak{g}^{e}(\geq1)$ is generated by
$\mathfrak{g}^{e}(1)$. Recall that the Lie superbracket $[\cdotp,\cdotp]:\mathfrak{g}_{\bar{1}}\times\mathfrak{g}_{\bar{1}}\rightarrow\mathfrak{g}_{\bar{0}}$
depends on $\sigma_{1}$, $\sigma_{2}$, $\sigma_{3}$ such that $\sigma_{1}+\sigma_{2}+\sigma_{3}=0$
and $\sigma_{i}\neq0$ for $i=1,2,3$. We compute that $[v_{1,1,-1}-v_{-1,1,1},v_{1,1,-1}-v_{-1,1,1}]=4\sigma_{2}E_{2}$,
$[v_{1,-1,1}-v_{-1,1,1},v_{1,-1,1}-v_{-1,1,1}]=4\sigma_{3}E_{3}$,
$[v_{1,1,-1}-v_{-1,1,1},v_{1,-1,1}-v_{-1,1,1}]=-2\sigma_{1}E_{1}+2\sigma_{2}E_{2}+2\sigma_{3}E_{3}$
and $[E_{2},v_{1,-1,1}-v_{-1,1,1}]=v_{1,1,1}$. Hence, we have that
$\mathfrak{g}^{e}(\geq1)$ is generated by $\mathfrak{g}^{e}(1)$. 

We have that $e=-\frac{1}{2\sigma_{1}}[v_{1,1,-1}-v_{-1,1,1},v_{1,-1,1}-v_{-1,1,1}]+\frac{\sigma_{1}+\sigma_{2}}{4\sigma_{1}\sigma_{2}}[v_{1,1,-1}-v_{-1,1,1},v_{1,1,-1}-v_{-1,1,1}]+\frac{\sigma_{1}+\sigma_{3}}{4\sigma_{1}\sigma_{3}}[v_{1,-1,1}-v_{-1,1,1},v_{1,-1,1}-v_{-1,1,1}]$.
Clearly $-\frac{1}{2\sigma_{1}}\neq0$ and since $\sigma_{i}\neq0$
for $i=1,2,3$, we have that $\frac{\sigma_{1}+\sigma_{2}}{4\sigma_{1}\sigma_{2}}=\frac{-\sigma_{3}}{4\sigma_{1}\sigma_{2}}\neq0$.
Similarly $\frac{\sigma_{1}+\sigma_{3}}{4\sigma_{1}\sigma_{3}}\neq0$.
Therefore, we deduce that $e\in[\mathfrak{g}^{e},\mathfrak{g}^{e}]$.
However, we have that $e$ is not strongly reachable as $\mathfrak{g}^{e}(1)\not\subseteq[\mathfrak{g}^{e},\mathfrak{g}^{e}]$.
\begin{itemize}
\item $\mathfrak{g}=G(3)$, $e=E+x_{1}$
\end{itemize}
According to Table \ref{tab:G(3)-2}, we know that $\mathfrak{g}^{e}=\mathfrak{g}^{e}(0)\oplus\mathfrak{g}^{e}(1)\oplus\mathfrak{g}^{e}(2)\oplus\mathfrak{g}^{e}(3)$
where basis elements of $\mathfrak{g}^{e}(j)$ for $j=0,1,2,3$ are
given in the following table. 
\begin{table}[H]
\noindent \begin{centering}
\begin{tabular}{|c|c|}
\hline 
$\mathfrak{g}^{e}(j)$ & Basis elements of $\mathfrak{g}^{e}(j)$\tabularnewline
\hline 
\hline 
$\mathfrak{g}^{e}(0)$ & $x_{6},y_{6},h_{1}+2h_{2},v_{1}\otimes e_{-3}-v_{-1}\otimes e_{-2},v_{1}\otimes e_{2}+v_{-1}\otimes e_{3}$\tabularnewline
\hline 
$\mathfrak{g}^{e}(1)$ & $v_{1}\otimes e_{0}-v_{-1}\otimes e_{1}$\tabularnewline
\hline 
$\mathfrak{g}^{e}(2)$ & $E,x_{1},v_{1}\otimes e_{3},v_{1}\otimes e_{-2}$\tabularnewline
\hline 
$\mathfrak{g}^{e}(3)$ & $x_{5},y_{2},v_{1}\otimes e_{1}$\tabularnewline
\hline 
\end{tabular}
\par\end{centering}
\caption{$\tau$-grading on $\mathfrak{g}^{e}$ when $e=E+x_{1}$}

\end{table}
 Since $\dim\mathfrak{g}^{e}(1)=1$, in order to check whether $e$
satisfies the Panyushev property, we only need to calculate $[v_{1}\otimes e_{0}-v_{-1}\otimes e_{1},v_{1}\otimes e_{0}-v_{-1}\otimes e_{1}]=-8E+8x_{1}$.
This implies that $\mathfrak{g}^{e}(\geq1)$ is not generated by $\mathfrak{g}^{e}(1)$
as we cannot obtain all basis elements of $\mathfrak{g}^{e}(2)$ from
commutators between basis elements of $\mathfrak{g}^{e}(1)$. Hence,
the Panyushev property does not hold for this case.

Note that $e\in[\mathfrak{g}^{e},\mathfrak{g}^{e}]$ if and only if
$e\in[\mathfrak{g}_{\bar{1}}^{e}(2),\mathfrak{g}_{\bar{1}}^{e}(0)]+[\mathfrak{g}_{\bar{1}}^{e}(1),\mathfrak{g}_{\bar{1}}^{e}(1)]$.
We have calculated that $[\mathfrak{g}_{\bar{1}}^{e}(1),\mathfrak{g}_{\bar{1}}^{e}(1)]=\langle-8E+8x_{1}\rangle$.
Then we calculate that $[v_{1}\otimes e_{3},v_{1}\otimes e_{-3}-v_{-1}\otimes e_{-2}]=[v_{1}\otimes e_{-2},v_{1}\otimes e_{2}+v_{-1}\otimes e_{3}]=16E+4x_{1}$.
This implies that $[\mathfrak{g}_{\bar{1}}^{e}(2),\mathfrak{g}_{\bar{1}}^{e}(0)]=\langle16E+4x_{1}\rangle$.
Thus we have that 
\[
e=\frac{3}{40}[v_{1}\otimes e_{0}-v_{-1}\otimes e_{1},v_{1}\otimes e_{0}-v_{-1}\otimes e_{1}]+\frac{1}{10}[v_{1}\otimes e_{3},v_{1}\otimes e_{-3}-v_{-1}\otimes e_{-2}].
\]
 The above equality holds if and only if $\mathrm{char}(\mathbb{K})=p\neq5$.
Hence, when $p\neq5$, we have that $e\in[\mathfrak{g}^{e},\mathfrak{g}^{e}]$,
i.e. $e$ is reachable. When $p=5$, we have that $e\notin[\mathfrak{g}^{e}(2),\mathfrak{g}^{e}(0)]+[\mathfrak{g}^{e}(1),\mathfrak{g}^{e}(1)]$
and thus $e$ is not reachable.

However, we have that $e$ is not strongly reachable because we cannot
obtain the basis element of $\mathfrak{g}^{e}(1)$ from $[\mathfrak{g}^{e},\mathfrak{g}^{e}]$. 
\begin{itemize}
\item $\mathfrak{g}=F(4)$, $e=E+e_{(3^{2},1)}$
\end{itemize}
Based on Tables \ref{tab:F(4)-centralizer} and \ref{tab:z(g^e)-F(4)},
we know that $\mathfrak{g}^{e}=\mathfrak{g}^{e}(0)\oplus\mathfrak{g}^{e}(1)\oplus\mathfrak{g}^{e}(2)\oplus\mathfrak{g}^{e}(3)\oplus\mathfrak{g}^{e}(4)$
where basis elements of $\mathfrak{g}^{e}(j)$ for $j=0,1,2,3,4$
are given in the following table. 

\begin{table}[H]
\noindent \begin{centering}
\begin{tabular}{|c|c|}
\hline 
$\mathfrak{g}^{e}(j)$ & Basis elements of $\mathfrak{g}^{e}(j)$\tabularnewline
\hline 
\hline 
$\mathfrak{g}^{e}(0)$ & $R_{1,-1}-R_{2,-2}+R_{3,-3}$\tabularnewline
\hline 
$\mathfrak{g}^{e}(1)$ & $v_{1}\otimes e_{1}s-v_{-1}\otimes e_{1}e_{2}e_{3}s,v_{1}\otimes e_{2}s,v_{-1}\otimes e_{1}e_{2}s+v_{1}\otimes e_{2}e_{3}s,v_{1}\otimes e_{1}e_{3}s$\tabularnewline
\hline 
$\mathfrak{g}^{e}(2)$ & $E,e_{(3^{2},1)},R_{2,-3},R_{2,0},R_{1,0},R_{1,3}$\tabularnewline
\hline 
$\mathfrak{g}^{e}(3)$ & $v_{1}\otimes e_{1}e_{2}s,v_{1}\otimes e_{1}e_{2}e_{3}s$\tabularnewline
\hline 
$\mathfrak{g}^{e}(4)$ & $R_{1,2}$\tabularnewline
\hline 
\end{tabular}
\par\end{centering}
\caption{$\tau$-grading on $\mathfrak{g}^{e}$ when $e=E+e_{(3^{2},1)}$}

\end{table}

Let $x=v_{1}\otimes e_{1}s-v_{-1}\otimes e_{1}e_{2}e_{3}s$ and $y=v_{-1}\otimes e_{1}e_{2}s+v_{1}\otimes e_{2}e_{3}s$.
We compute $[x,x]=R_{1,0}$, $[x,v_{1}\otimes e_{2}s]=\frac{1}{2}R_{2,0}$,
$[x,v_{1}\otimes e_{1}e_{3}s]=R_{1,3}$, $[x,y]=e_{(3^{2},1)}-6E$,
$[v_{1}\otimes e_{2}s,v_{1}\otimes e_{1}e_{3}s]=6E$ and $[v_{1}\otimes e_{2}s,y]=R_{2,-3}$.
Thus we have that $[\mathfrak{g}^{e}(1),\mathfrak{g}^{e}(1)]=\mathfrak{g}^{e}(2)$.
We further calculate that $[R_{1,0},v_{1}\otimes e_{2}s]=-v_{1}\otimes e_{1}e_{2}s$
and $[R_{1,0},y]=v_{1}\otimes e_{1}e_{2}e_{3}s$. Hence, we have that
$[\mathfrak{g}^{e}(2),\mathfrak{g}^{e}(1)]=\mathfrak{g}^{e}(3)$.
Similarly $[\mathfrak{g}^{e}(3),\mathfrak{g}^{e}(1)]=\mathfrak{g}^{e}(4)$
as $[v_{1}\otimes e_{1}e_{2}e_{3}s,y]=R_{1,2}$. Therefore, we conclude
that $\mathfrak{g}^{e}(\geq1)$ is generated by $\mathfrak{g}^{e}(1)$. 

Note that $e=[x,y]+\frac{7}{6}[v_{1}\otimes e_{2}s,v_{1}\otimes e_{1}e_{3}s]$.
Hence, we have that $e$ is reachable. 

Since $\mathfrak{g}^{e}(0)$ is $1$-dimensional, we have that the
subspace of $[\mathfrak{g}^{e},\mathfrak{g}^{e}]$ with grading $0$
is an empty set. Hence, we cannot obtain the basis element of $\mathfrak{g}^{e}(0)$
from $[\mathfrak{g}^{e},\mathfrak{g}^{e}]$ and thus $e$ is not strongly
reachable.

\subsection*{Statements and Declarations. }

The author declares that she has no conflict of interest.

\end{document}